\documentclass[final,onefignum,onetabnum]{siamart251216}

\usepackage{lipsum}
\usepackage{amsfonts}
\usepackage{graphicx}
\usepackage{epstopdf}
\usepackage{algorithmic}
\ifpdf
  \DeclareGraphicsExtensions{.eps,.pdf,.png,.jpg}
\else
  \DeclareGraphicsExtensions{.eps}
\fi

\usepackage{enumitem}
\setlist[enumerate]{leftmargin=.5in}
\setlist[itemize]{leftmargin=.5in}

\newsiamremark{remark}{Remark}
\newsiamremark{hypothesis}{Hypothesis}
\crefname{hypothesis}{Hypothesis}{Hypotheses}
\newsiamthm{claim}{Claim}
\newsiamremark{fact}{Fact}
\crefname{fact}{Fact}{Facts}

\headers{An Example Article}{D. Doe, P. T. Frank, and J. E. Smith}

\title{An Example Article\thanks{Submitted to the editors DATE.
\funding{This work was funded by the Fog Research Institute under contract no.~FRI-454.}}}

\author{Dianne Doe\thanks{Imagination Corp., Chicago, IL 
  (\email{ddoe@imag.com}, \url{http://www.imag.com/\string~ddoe/}).}
\and Paul T. Frank\thanks{Department of Applied Mathematics, Fictional University, Boise, ID 
  (\email{ptfrank@fictional.edu}, \email{jesmith@fictional.edu}).}
\and Jane E. Smith\footnotemark[3]}

\usepackage{amsopn}

\ifpdf
\hypersetup{
  pdftitle={Higher-Order Moments of Correlated Multivariate Hermite Polynomials},
  pdfauthor={L.~Lyman}
}
\fi

\usepackage{graphicx,epstopdf}
\usepackage{amssymb,url}
\usepackage{mathtools}      
\usepackage{nccmath}        
\usepackage{bm}
\usepackage{array}
\usepackage{booktabs}       
\usepackage{float}
\usepackage{euscript}       
\usepackage{xcolor}

\newcommand{\N}{\mathbb{N}}
\newcommand{\R}{\mathbb{R}}

\newcommand{\balpha}{\bm{\alpha}}

\newcommand{\oalpha}{|\balpha|}
\newcommand{\bbeta}{\bm{\beta}}
\newcommand{\obeta}{|\bm{\beta|}}
\newcommand{\bgamma}{\bm{\gamma}}
\newcommand{\ogamma}{|\bgamma|}
\newcommand{\bdelta}{\bm{\delta}}
\newcommand{\odelta}{|\bdelta|}
\newcommand{\A}{\bm{A}}
\newcommand{\B}{\bm{B}}

\newcommand{\x}{\bm{x}}
\newcommand{\bxi}{\bm{\xi}}
\newcommand{\doub}{\langle H_{\balpha}, H_{\bbeta} \rangle}
\newcommand{\trip}{\langle H_{\balpha}, H_{\bbeta}, H_{\bgamma} \rangle}

\newsiamthm{example}{Example}

\allowdisplaybreaks

\title{Higher-Order Moments of Correlated Multivariate Hermite Polynomials%
\thanks{\funding{Clare Booth Luce Fellowship (2024--present).}}}

\author{Laura Lyman\thanks{Department of Mathematics and Statistics, Mount Holyoke College (\email{lalyman@mtholyoke.edu}).}}

\begin{document}

\headers{Higher-Order Moments of Correlated Hermite Polynomials}{Laura Lyman}

\maketitle

\begin{abstract} Stochastic Galerkin methods reformulate a differential equation with uncertain parameters as a deterministic surrogate system by writing the solution as a series of polynomial basis functionals and projecting onto the parameter space. Given correlated Gaussian parameters, Rahman (2017) established solution convergence when the basis functionals are multivariate Hermite polynomials (MHPs), which is the setting considered here. Each term of the original equation contributes coefficients to the surrogate system, which we call the $N$th moments, that are expectations of products of the basis polynomials. Since the parameters are dependent, the $N$th moments do not factor into one-dimensional expectations via tensorization as they do in the independent case; consequently, their values are typically determined numerically at a nontrivial and compounding cost. Namely, the number of $N$th moments that emerge grows combinatorially with parameter dimension and truncation order. In this paper, we first relate MHPs and slowly growing functions --- which are smooth functions whose derivatives of all orders grow at most polynomially --- in correlated Gaussian parameters. In doing so, we derive an expression for the $N$th moment of an MHPs in terms of $(N-1)$th moments that can be evaluated exactly, along with sufficient conditions for the $N$th moment to vanish. Numerical examples illustrate the benefit of evaluating such moments directly without numerical approximation and confirm that formula-predicted values match sampling-based estimates.

\end{abstract}




\section{Introduction}\label{sec:intro}
Differential equations that model natural phenomena often contain terms that are only known up to some level of certainty. Within the framework of \emph{forward propagation problems}, this uncertainty is typically represented by a random vector with prescribed joint density that is incorporated into the model equation, whose solution then becomes stochastic. A basic objective of uncertainty quantification (UQ) is to describe the distribution of the solution at an affordable cost \cite{UQHandbook_IntroUQ, Smith}. Stochastic Galerkin schemes (SGS) are especially useful when traditional sampling-based approaches are prohibitively expensive to implement, which can occur if resolving the model equation for a single realization of the random parameter is costly. In particular, SGS express the solution as a series of polynomial functionals of the random vector, whose deterministic expansion coefficients (the \emph{random modes}) are themselves the solution of a surrogate system of deterministic differential equations  \cite{alexanderian2013:master, constantine_thesis:master, ghanem_and_spanos:master, le_maitre:master, GI_book}. For a thorough introduction of SGS and their context within UQ, see \cite{mytripprodpaper}.

The basis functions in the SGS expansion are customarily polynomials that are orthogonal with respect to the joint density of the random vector, in which case the series is often called a \emph{polynomial chaos expansion} (PCE).  Wiener \cite{wiener:master} first wrote a Gaussian process as an infinite series of Hermite polynomials, which Ghanem and Spanos \cite{ghanem_and_spanos:master} later truncated for the sake of approximation.  This mathematical setting is reviewed in \cref{SGS_intro}.

When compared to MC algorithms, initial spectral method applications to UQ problems showed \emph{orders-of-magnitude} reductions in the cost needed to estimate statistics with comparable accuracy \cite{xiu_karniadakis:master}; however, there are of course some exceptions, e.g., see \cite{exceptionsToPCefficiency1, exceptionsToPCefficiency2}. Several PCE extensions have been developed in the last years  \cite{Navarro} that have shown to be effective in a large number of applications in different fields: fluid dynamics \cite{Feinberg,le_maitre:master, verif_valid_ref_2, GI_book}, air-water flows in soil \cite{airflowsoil}, aircraft design \cite{ZDR}, radiation therapy \cite{Perko_2016, proton_therapy_thesis}, etc. Further, spectral Galerkin methods have the advantage of returning the entire distribution of the output QoI, from which any descriptive statistics (mean, variance, etc.) can be computed directly; this is in direct contrast to sampling strategies, which only compute pre-specified descriptive statistics.


The price of these benefits is that SGS are \emph{intrusive}: the projected system must be derived and implemented from scratch, and assembling it requires expectations of products of basis polynomials that never appear in a deterministic code \cite{Debusschere2017, GI_book, Smith}. (When an existing solver must be treated as a black box, non-intrusive alternatives such as MC sampling or stochastic collocation \cite{xiu_stochcollsurvey} are usually preferred \cite{constantine_thesis:master}.) The products in question are dictated by the nonlinearity of the model. A linear term leads only to inner products of two basis polynomials, but a term of polynomial degree $p$ in the unknowns leads, after projection, to expectations of products of $p+1$ basis polynomials. We call such an expectation of $N$ polynomials an \emph{$N$th moment} or \emph{$N$th product}. When the components of $\bxi$ are independent, every basis polynomial is a product of univariate orthogonal polynomials, so each $N$th moment splits into one-dimensional integrals; these vanish by orthogonality or are given by classical formulas, such as the triple product formula \cref{trip_prod_1d} for Hermite polynomials.

 In the multivariate case, one of the \emph{essential} assumptions when applying SGS is that the stochastic model parameters are independent. This assumption allows the joint distribution of the random variables to decompose as a product of probability densities, enabling the easy construction of multivariate orthogonal polynomials via tensorizing the spaces of univariate orthogonal polynomials \cite{rahman:master, Smith, xiu_thesis:master}. These orthogonal basis functionals are consistent with the underlying measure \cite{ernst:master} and thus constitute a complete Hilbert basis \cite{rahman:master}; as clarified only recently \cite{ernst:master}, this completeness\footnote{The completeness of the Hilbert basis is equivalent to the underlying probability measure being uniquely determined by its moments \cite{ernst:master}.} is vitally important for convergence of the truncated PCE to the correct limit. Clasically, decorrelating transformations \cite{Rosenblatt, Nataf} introduce numerical instabilities --- and while this is an active area of research, this paper does not consider transformation-based approaches.
 
 In the  special case of correlated Gaussian inputs, $\bxi \sim \mathcal{N}(\bm{0}, \Sigma)$, Rahman \cite{rahman:master} showed that the \emph{multivariate Hermite polynomials} (MHPs) associated with the joint density form a complete basis, and that the resulting \emph{correlated} polynomial chaos expansion (CPCE) converges in mean square, in probability, and in distribution. The catch is that MHPs are only \emph{weakly orthogonal} \cite{dunkl_xu_2014:master, rahman:master}: polynomials of different total degree are orthogonal, but polynomials of the same degree generally are not. Since the joint density no longer factors, neither do the moments. Even second moments call for a dedicated formula \cite{mydoubprodpaper}, and prior to this work no closed-form expression was known for $N \ge 3$. The usual recourse is brute-force sampling, which scales poorly. For the delayed logistic model of \cref{example_application}, a CPCE with only $n = 2$ random inputs and total degree $M = 2$ already involves up to $189$ double, triple, and quadruple products, and $n = 3$ raises the count to $560$. Namely, the number of $N$th moments that emerge during Galerkin projection when the parameters are correlated grows combinatorially with parameter dimension and truncation order.

\subsection{Summary of Contributions}
This paper develops the theory needed to evaluate $N$th moments of MHPs in correlated Gaussian variables without numerical approximation. Throughout the document, $\bxi \sim \mathcal{N}(\bm{0}, \Sigma)$ with a known symmetric positive definite (SPD) covariance matrix $\Sigma$. Every expectation involving MHPs is taken with respect to the law of the $\bxi$ random vector. Our main results are the following.

\begin{itemize}
\item The expectation of any slowly growing function (\Cref{slowly_increasing}) of $\bxi$ --- that is, any smooth function whose derivatives of all orders grow at most polynomially  --- is finite. In particular, every $N$th moment of MHPs exists (\Cref{moments-exist}).
\item For every multi-index $\A$ and every slowly growing $f$, we have $\mathbb{E}(H_{\A}(\bxi) f(\bxi)) = \mathbb{E}(\partial^{\A} f(\bxi))$ (\cref{IBP}). Consequently, any one factor of an $N$th product may be traded for a derivative acting on the remaining $N-1$ factors (\cref{general_deriv_N_product}).
\item There is an explicit, recursive formula for the product of $N$ multivariate Hermite polynomials in correlated Gaussian random variables; namely, for $\oalpha \le \obeta \le |\bgamma^{(1)}| \le \cdots \le |\bgamma^{(N-2)}|$ (\cref{N_prod}),
\begin{multline*}
\langle H_{\balpha}, H_{\bbeta}, H_{\bgamma^{(1)}}, \ldots, H_{\bgamma^{(N-2)}} \rangle \\
= \sum_{\A \le \balpha} \sum_{\substack{\B \le \bbeta \\ |\B| = |\A|}} \mbinom{\balpha}{\A} \mbinom{\bbeta}{\B} \langle H_{\A}, H_{\B} \rangle \langle H_{\balpha + \bbeta - \A - \B}, H_{\bgamma^{(1)}}, \ldots, H_{\bgamma^{(N-2)}} \rangle.
\end{multline*}
Further, this expression reduces to the triple product ($N = 3$) formula recently provided \cite{mytripprodpaper}. 
\item The $N$th product is guaranteed to be zero whenever the order of its largest index exceeds the sum of the orders of the others (\cref{N_prod_zeros}).
\end{itemize}
\medskip

This paper is organized as follows. \Cref{math_background} presents relevant definitions, notational conventions, and acronyms used throughout the document. \Cref{SGS_math} describes the formal probabilistic setting. In \cref{SGS_intro}, we overview spectral Galerkin expansions and the original series polynomial chaos expansion introduced by Wiener.  \Cref{MHP_section} introduces correlated MHPs, which serve in \cref{sec_PC_corr_gauss} as the basis of the CPCE.  \Cref{example_application} presents an example application of solving a stochastic
variant of Hutchinson's equation for population growth  with the CPCE and helps illustrate the practical value
of a triple product formula that a computer can evaluate directly. Then in \cref{Nthprodformula} we present the paper’s main result: an $N$th product formula that can be directly evaluated, after which we work through some specific example calculations.  \Cref{simple_MC_conv} compares the formula with MC estimates, and \cref{sec:conclusions} concludes. The existence results and various identities are proved in \cref{appendix_defs,existence_section}.

\section{Definitions and notation}\label{math_background}
We write $\N = \{1, 2, \ldots\}$, $\N_0 = \{0, 1, 2, \ldots\}$, and $[n] = \{1, \ldots, n\}$ for $n \in \N$. The acronyms listed in \cref{tab_acronyms} are used throughout the document. 

\begin{table}[htbp]
\caption{Acronyms used throughoug the paper.}\label{tab_acronyms}
\centering
\begin{tabular}{ll}
\toprule
SGS & stochastic Galerkin scheme \\
PCE & polynomial chaos expansion \\
MHP & multivariate Hermite polynomial \\
CPCE & correlated polynomial chaos expansion \\
\bottomrule
\end{tabular}
\end{table}
We now present notational conventions and definitions similar to those in \cite{mytripprodpaper}. 
\begin{definition}[Multi-index]\label{multi-index} An $n$-dimensional \textit{multi-index} $\balpha$ is an $n$-tuple $(\alpha_1, \ldots, \alpha_n) \in \N_0^n$ of non-negative integers.  Each $\alpha_i$ is referred to as the \textit{$i$th element of the multi-index~$\balpha$.}
\end{definition}

\noindent Unless otherwise specified, a multi-index $\balpha$ is assumed to have the notation given in Def.~\ref{multi-index}. We also write $\bm{e}_j$ for the $j$th unit multi-index whose $j$th entry is $1$ and whose other entries are $0$.

For multi-indices $\balpha, \bbeta \in \N_0^n$, one defines \cite{rahman:master,slepian:master}
\begin{enumerate}
\item \emph{componentwise sum and difference}: $\balpha \pm \bbeta = (\alpha_1 \pm \beta_1, \ldots, \alpha_n \pm \beta_n)$;
\item \emph{absolute value}: $|\balpha| := \|\balpha\|_1 = \alpha_1 + \cdots + \alpha_n$, which we call the \emph{order} of $\balpha$;
\item \emph{factorial}: $\balpha ! = \alpha_1! \cdots \alpha_n! = \prod_{i=1}^n \alpha_i !$;
\item \emph{partial derivative}: $\partial^{\balpha} = \partial_1^{\alpha_1} \cdots \partial_n^{\alpha_n}$.
\end{enumerate}

Moving forward, any $\balpha$ will denote a multi-index of order $k \ge 1$ unless otherwise specified. Following  \cite{rahman:master}, we use the convention of \emph{graded lexicographic ordering}  throughout this document. 

\begin{definition}[Graded lexicographic ordering]\label{glex}For multi-indices $\balpha \in \N_0^{n}$, we use a \emph{graded lexicographic ordering}; that is
\[ \balpha \le \bbeta \Leftrightarrow  |\balpha| < |\bbeta| \text{ or } \left( |\balpha| = |\bbeta| \text{ and the \emph{leftmost} digit of } \balpha - \bbeta \text{ is positive}\right). \] 
\end{definition}
For instance, for $n = 2$,
\vspace{-3pt}
\[ (0,0) < \overbrace{(1,0) < (0,1)}^{|\balpha| = 1} < \overbrace{(2,0) < (1,1) < (0,2)}^{|\balpha| = 2} <  \overbrace{(3,0) < (2,1) < (1,2) < (0, 3)}^{|\balpha| = 3} < \cdots\] 
In all matrices, vectors, and sums indexed by $\balpha$, we assume a graded lexicographic ordering (in ascending order) when needed. 

For multi-indices, Leibniz' formula \cite{leibniz-rule:master} for the generalized product rule is as follows:
\begin{equation}\label{leibniz-rule}
\partial^{\balpha}(fg) = \sum_{\A \le \balpha} \mbinom{\balpha}{\A} \partial^{\A} (f) \partial^{\balpha - \A} (g).
\end{equation}

With this notation in place, we review the mathematical theory behind multivariate spectral expansions.

\section{Probabilistic setting}\label{SGS_math}

\subsection{Input}\label{input}
Let $(\Omega, \mathcal{A}, \mathbb{P})$ be a complete probability space, where $\Omega$ is a nonempty sample space of elementary events, $\mathcal{A}$ is a Borel $\sigma$-algebra on $\Omega$, and $\mathbb{P}: \mathcal{A} \rightarrow [0,1]$ is a probability measure. For some $\mathbb{A}^n \subseteq \R^n$, suppose we are given an $\mathbb{A}^n$-valued input random vector  $\bm{\xi}: (\Omega, \mathcal{A}) \rightarrow (\mathbb{A}^n, \mathcal{B}_{\mathbb{A}^n})$ and its density function $\varphi: \mathbb{A}^n \rightarrow \R_{\ge 0}$, where $\mathcal{B}_{\mathbb{A}^n}$ denotes the Borel sigma algebra on $\mathbb{A}^n$. The positive integer $n$ represents the number of input random variables (i.e., the number of components $\xi_1, \ldots, \xi_n$ of $\bxi$) and is sometimes referred to as the \emph{dimension} of the stochastic problem \cite{RAHMAN2019124538}.

\subsection{Output and $N$th product definition}\label{output}

Given the random vector $\bm{\xi}$ and its density $\varphi$ on $\mathbb{A}^n$, consider some measurable function $f:(\mathbb{A}^n, \mathcal{B}_{\mathbb{A}^n}) \rightarrow (\R, \mathcal{B}_{\R})$, where $\mathcal{B}_{\R}$ is the Borel sigma algebra on the real line. Then per usual, $f(\bm{\xi}): (\Omega, \mathcal{A}) \rightarrow (\R, \mathcal{B}_{\R})$ is the random variable defined by $f(\bm{\xi})(\omega) = f(\bm{\xi}(\omega)).$  

Consider the Hilbert space\footnote{The precise definition of $L^p (\Omega, \mathcal{A}, \mathbb{P})$ for $1 \le p < \infty$ requires introducing an equivalence relation, in which $f$ and $g$ are considered equivalent if $f = g$ almost everywhere. Therefore, $L^2(\Omega, \mathcal{A}, \mathbb{P})$ technically consists of \emph{equivalence classes} of $(\Omega, \mathcal{A})$-measurable functions that have finite second moments. However, we follow the standard practice of describing elements of $L^2$ as functions rather than equivalence classes of functions, and all statements/equalities are understood to hold almost everywhere (a.e.)~\cite{ernst:master,Lp_def:master}.} of square-integrable functions,

\begin{align*}
\hspace{-5pt} L^2(\Omega, \mathcal{A}, \mathbb{P}) = \left\{ f: \Omega \rightarrow \R \;\middle|\; \hspace{-2pt} f \text{ is } (\Omega, \mathcal{A})\text{-measurable and} \hspace{-1pt}  \int_\Omega |f(\omega)|^2 \, d \mathbb{P}(\omega) < \infty \right\}\hspace{-2pt}.
\end{align*}
Following common convention, we exploit the equivalence between realizations $\omega \in \Omega$ and the output values of $\bm{\xi}(\omega) \in \mathbb{A}^n$ to pose problems in $(\Omega, \mathcal{A})$ instead in the image probability space $(\mathbb{A}^n, \mathcal{B}_{\mathbb{A}^n}, \varphi(\bm{y}) d\bm{y})$ \cite{rahman:master}. Accordingly, it can be shown\footnote{This follows from the Doob-Dynkin lemma, as described in \cite{ernst:master}.} that $f(\bm{\xi}) \in L^2(\Omega, \mathcal{A}, \mathbb{P})$ if and only if $f \in L^2(\mathbb{A}^n, \mathcal{B}_{\mathbb{A}^n}, \varphi(y) dy)$, where
\begin{align*}
\begin{split}
&L^2(\mathbb{A}^n, \mathcal{B}_{\mathbb{A}^n}, \varphi(\bm{y})\, d \bm{y}) \\
&\quad = \left\{ f: \mathbb{A}^n \rightarrow \R \;\middle|\; f \text{ is } (\mathbb{A}^n, \mathcal{B}_{\mathbb{A}^n})\text{-measurable and } \int_{\mathbb{A}^n} |f(\bm{y})|^2 \varphi(\bm{y})\, d\bm{y}< \infty \right\}.
\end{split}
\end{align*}

 \noindent Therefore, we can examine the convergence of random variables $f(\bxi) \in L^2(\Omega, \mathcal{A}, \mathbb{P})$ and functions $f \in L^2(\mathbb{A}^n, B_{\mathbb{A}^n}, \varphi(\bm{y}) d \bm{y})$ interchangeably.

One can construct an inner product on either of these $L^2$ spaces:

\begin{align}
\langle f(\bm{\xi}), g(\bm{\xi}) \rangle_{d\mathbb{P}(\omega)}&:= \int_{\Omega} f(\bxi(\omega)) g(\bxi(\omega))  d \mathbb{P}(\omega) \label{inner_prod}  \\
&\phantom{:}= \int_{\mathbb{A}^n} f(\bm{y}) g(\bm{y}) \varphi(\bm{y}) d\bm{y} =: \langle f, g \rangle_{\varphi(\bm{y}) d\bm{y}} \notag
\end{align}
where \cref{inner_prod} also equals  $\mathbb{E}_{\bm{\xi}}(f(\bxi) g(\bxi))$. The inner product induces the norm $\|f \|_{\varphi(\bm{y}) d\bm{y}}^2 := \langle f^2 \rangle_{\varphi(\bm{y}) d\bm{y}}$ \cite{constantine_primer:master, Amir:master}. We often omit the subscript on the angled inner product brackets. Sometimes, we refer to $\langle f, g \rangle$ as a \emph{double product}. The \emph{triple product} $\langle f, g, h \rangle$, \emph{quadruple product} $\langle f_1, f_2, f_3, f_4 \rangle$, and $N$th  product $\langle f_1, \ldots, f_N \rangle$ are defined analogously.

\section{Spectral Galerkin expansions}\label{SGS_intro}

 Consider solving a differential equation that incorporates the $\bm{\xi}$ stochastic quantity. That is, 

\begin{equation}\label{du_dt}
 \mathcal{F}(u, \bm{d}; \bm{\xi}) =  0, \qquad \bm{d} \in \mathcal{D} \subset \R^p.
 \end{equation}
The $\mathcal{D}$ is a nonempty, compact domain representing the space of the $p$ deterministic equation parameters. Here, $\mathcal{F}$ is a differential operator, often containing linear and nonlinear terms. The underlying assumption is that $\mathcal{F}$ is prohibitively expensive to evaluate directly with sampling-intensive methods, so it is useful to identify a less expensive model whose input-output behavior closely mimics that of the underlying forward function \cite{mueller2024polynomialchaossurrogateconstruction}. Consequently, the reduced-cost model acts as a \emph{surrogate representation} for $\mathcal{F}$ that can be used instead to study uncertainty propagation.

In (\ref{du_dt}), we have that $u:\mathcal{D} \times S_{\bm{\xi}} \rightarrow \R$, where $S_{\bm{\xi}} := \text{supp}(\bm{\xi})$ is the \emph{support} of the $\bm{\xi}$ random vector. For each $\bm{d} \in \mathcal{D}$ fixed, $u(\bm{d}; \bm{\xi}): (\Omega, \mathcal{A}) \rightarrow (\R, \mathcal{B}_{\R}) $ denotes the random variable defined by $u(\bm{d};\bm{\xi})(\omega) = u(\bm{d}, \bm{\xi}(\omega))$. Then $\{ u(\bm{d}; \bm{\xi}) \}_{\bm{d} \in \mathcal{D}}$  is a \emph{stochastic process}  in the sense that it is a collection of random variables indexed by $\bm{d} \in \mathcal{D}$ \cite{Amir:master, Smith}. 
\emph{Following a common abuse of notation, we may write $u(\bm{d};\bm{\xi})$ to mean $\{ u(\bm{d}; \bm{\xi}) \}_{\bm{d} \in \mathcal{D}}$ when $\bm{d} \in \mathcal{D}$ is unspecified. Similarly, we will often use the terms $u, u(\bm{d};\bm{\xi}),$ and  $\{ u(\bm{d}; \bm{\xi}) \}_{\bm{d} \in \mathcal{D}}$ interchangeably when referring to the ``solution" of (\ref{du_dt}).}

A \emph{spectral Galerkin representation}, sometimes called a just a \emph{spectral expansion}, \cite{constantine_thesis:master, le_maitre:master, stoch_galer_2:master, GI_book, salazar, Smith, Soize_Ghanem, xiu_thesis:master, xiu_karniadakis:master} of a stochastic process $u(\bm{d}; \bm{\xi})$ aims to find a series expansion of the form
\begin{equation}\label{high_dim_PCE}
u(\bm{d}; \bm{\xi}) = \sum_{\balpha \in \N_0^n} \thinspace u_{\balpha}(\bm{d}) \Psi_{\balpha}(\bm{\xi})
\end{equation}where the $\balpha = (\alpha_1, \ldots, \alpha_n)$ are multi-indices (Definition \ref{multi-index}) assumed in graded lexicographic order (Definition \ref{glex}).
In (\ref{high_dim_PCE}), 
\begin{itemize}
\item the $u_{\balpha}: \mathcal{D} \rightarrow \R$ are \emph{deterministic coefficient functions}, and 
\item the $\Psi_{\balpha}(\bm{\xi})$ are basis functionals.
\end{itemize}
The coefficient functions $u_i$ are also called \emph{random nodes} \cite{xiu_thesis:master, xiu_karniadakis:master}. By utilizing a spectral expansion, the randomness is effectively transferred onto the basis functionals. Namely, we can observe how the deterministic $\bm{d}$ and stochastic $\xi$ parameters of $u(\bm{d}; \xi)$ are separated within the series of (\ref{high_dim_PCE}). 

Typically, the functionals $\Psi_{\balpha}$ in (\ref{high_dim_PCE}) are chosen to be polynomials that are \emph{orthogonal} with respect to the $\bxi$ density --- that is, orthogonal with respect to the inner product defined in (\ref{inner_prod}) --- so that they constitute a \emph{complete} Hilbert basis; as clarified in \cite{ernst:master}, this completeness is vitally important for the convergence of the spectral representation. We will discuss this more in section \ref{trans_free_approaches}. 


\subsection{Wiener--Hermite polynomial chaos} Norbert Wiener first represented a Gaussian process as a countably infinite series of Hermite polynomials \cite{wiener:master}. The \emph{probabilist's Hermite polynomials} $H_i \in C^{\infty}(\R)$ are defined via the formula \cite{Digression_Herm_Polys}
\begin{align}\label{herm_1D_def}
H_i(x) = (-1)^i e^{x^2/2} \partial_x^i \left( e^{-x^2/2} \right).
\end{align}
The first few Hermite polynomials are given in Table \ref{herm_table}. 

We follow the notation in sections \ref{input} and \ref{output} and consider the special case of when $\xi$ is univariate standard normal (where ``$\bxi$" is written as ``$\xi$" since $n = 1$). Equation \ref{herm_1D_def} indicates that
\begin{align}\label{1d_herm_density}
H_i(\xi) = (-1)^i \varphi^{-1}(\xi) \partial_\xi^i \left( \varphi(\xi) \right),
\end{align}
where $\varphi:(\R, \mathcal{B}_\R) \rightarrow \R_{\ge 0}$ is the density of $\xi \sim \mathcal{N}(0,1)$. The statistics $H_i(\xi)$ are known \cite{constantine_primer:master, GI_book, Smith, polys:master} to have the following first order moments:
\begin{align}
\langle H_i \rangle &= \mathbb{E}_{\xi}(H_i(\xi)) = \delta_{i 0}, \label{single_prod_1D}
\end{align}
second order moments (\emph{double products}):
\begin{align}
 \langle H_i, H_j \rangle = \mathbb{E}_{\xi} (H_i(\xi) H_j(\xi) ) = i! \delta_{ij}, \label{double_prod_1D} 
\end{align}
and  third order moments (\emph{triple products}): 

\begin{align}
\langle H_i, H_j, H_k \rangle &= \mathbb{E}_{\xi}(H_i(\xi) H_j(\xi) H_k(\xi)) \notag \\
&= \begin{cases} 0 & \text{ if } i + j + k \text{ is odd or } \max(i,j,k) > s \label{trip_prod_1d} \\
\frac{i! j! k!}{(s-i)! (s-j)! (s-k)!} &\text{ else.} \end{cases} 
\end{align}
In (\ref{single_prod_1D}) and (\ref{double_prod_1D}), $\delta_{ij}$ is the Kronecker delta (which equals 1 for $i = j$ and 0 otherwise). In (\ref{trip_prod_1d}), $s := ( i+j+k)/{2}$ for $i,j,k \in \N_0$.

The Hermite polynomials $\{H_i(\xi)\}_{i=0}^\infty$ form a complete orthogonal basis of $L^2(\Omega, \mathcal{A}, \mathbb{P})$ \cite{ernst:master, rahman:master}. The double product formula \cref{double_prod_1D} is precisely this orthogonality property with respect to the standard normal density.


\begin{table}[htbp]
\caption{The probabilist's Hermite polynomials $H_i$ for $i = 0, \ldots, 5$.}\label{herm_table}
\centering
\renewcommand{\arraystretch}{1.15}
\begin{tabular}{cl}
\toprule
$i$ & $H_i(x)$ \\
\midrule
0 & $1$ \\
1 & $x$ \\
2 & $x^2 - 1$ \\
3 & $x^3 - 3x$ \\
4 & $x^4 - 6x^2 + 3$ \\
5 & $x^5 - 10 x^3 + 15 x$ \\
\bottomrule
\end{tabular}
\end{table}
The spectral expansion first proposed by Wiener,
\begin{align}\label{PC}
u(\bm{d}; \xi) = \sum_{i=0}^\infty u_i(\bm{d}) H_i(\xi),
\end{align}
is called \emph{Wiener-Askey chaos}, a \emph{Wiener-Hermite series}, a \emph{polynomial chaos expansion} (PCE), or simply a \emph{chaos expansion} \cite{mypaper:master, GI_book, rahman:master, Smith, xiu_thesis:master}. In this document, the \emph{PC expansion} or the \emph{(classic) PCE} will refer to the representation given by (\ref{PC}).

 Define the truncated PCE as
\begin{align}\label{S_M}
S_M(u(\bm{d}; \xi)) := \sum_{k=0}^M u_k(\bm{d}) H_k(\xi),
\end{align}
which is the projection of $u$ onto the space $\mathcal{P}_{M}(\xi)$ of all polynomials in $\xi$ of degree up to and including $M$ \cite{Smith}. If $\{u(\bm{d};\xi)\}_{d \in \mathcal{D}}$ is a \emph{second-order stochastic process}, i.e., $u(\bm{d};\xi) \in L^2(\Omega, \mathcal{A}, \mathbb{P})$ for all $\bm{d} \in \mathcal{D}$, then the Cameron-Martin theorem\footnote{The original argument by Cameron and Martin \cite{cammartin:master} is in the context of $L^2(\mathcal{C})$, where $\mathcal{C} = \{ x \in C[0,1] \mid x(0) = 0\}$; however, the theorem is later proved generally on the Hilbert space $L^2(\Omega, \mathcal{A}, \mathbb{P})$ \cite{ernst:master, janson_1997:master}.} states that the polynomial chaos expansion converges in mean square;\footnote{Note that Equation \ref{S_M_convergence} does  not say directly that $S_M(u(\bm{d}; \xi))$ as a stochastic \emph{process} converges to the stochastic process $u(\bm{d}; \xi)$; rather, the mean-square convergence refers to $S_M$ as a \emph{random variable} at some $\bm{d} \in \mathcal{D}$ converging to the random variable $u(\bm{d}; \xi)$ at the same fixed but arbitrary $\bm{d}$ index.} that is, for each fixed $\bm{d} \in \mathcal{D}$,
\begin{equation}\label{S_M_convergence}
S_M(u(\bm{d}; \xi)) \underset{L^2}{\overset{M \rightarrow \infty} {\xrightarrow{\hspace*{.8cm}}}} u(\bm{d};\xi).
\end{equation}
This justifies the PCE and its truncation to a finite number of terms for the sake of computation, as was first considered by Ghanem and Spanos \cite{ghanem_and_spanos:master} in a stochastic finite element method.

 \subsection{Transformation-free approaches}\label{trans_free_approaches}
There have been several transformation-free approaches to account for arbitrary but dependent Gaussian random variables. Using Gram--Schmidt orthogonalization \cite{Golub}, Navarro et al.~discussed construction of a PCE with multivariate orthogonal polynomials for correlated variables \cite{Navarro}. However, the existence of multivariate Hermite orthogonal polynomials, which can be used as a basis for dependent Gaussian measures, has not been recognized \cite{rahman:master}. Soize and Ghanem proposed orthogonal functional bases with regard to a general dependent probability measure of random input, but the bases elements are not necessarily \emph{polynomials} \cite{Soize_Ghanem}; as a result, an analytical treatment of PCE-based statistics in this case is highly non-trivial \cite{Smith}. Furthermore, both the works of Navarro and Soize/Ghanem do not address a \emph{fundamental} concern raised by Ernst et al.~about the completeness of their orthogonal basis functionals, which was shown to be a necessary condition for convergence \cite{ernst:master, rahman:master}. 

In 2017, Rahman constructed a PCE with \emph{weakly} orthogonal (defined in Eq.~\ref{weak_orthog}) basis functionals consistent with the Gaussian measure, which was shown to converge to the correct limit in mean-square, probability, and distribution  \cite{rahman:master}. We will refer to this series expansion as the \emph{correlated multivariate PCE} (or CPCE)\footnote{In \cite{rahman:master}, Rahman refers to the introduced countably infinite PC series in correlated Gaussians as the \emph{generalized PCE}. However, this language might cause confusion, given that the phrase \emph{generalized PCE} (gPCE) typically refers to the univariate PC expansion first introduced by Xiu and Karniadakis in 2002 for non-Gaussian random variables \cite{xiu_thesis:master, xiu_karniadakis:master}; see the discussion in \cite{mytripprodpaper}.} and use it exclusively throughout the rest of this document.

\section{Correlated multivariate Hermite polynomials}\label{MHP_section}
As in \cref{input}, let $(\Omega, \mathcal{A}, \mathbb{P})$ be a complete probability space. Let $\bm{\xi} \sim \mathcal{N}(\bm{0}, \Sigma)$ be an $n$-dimensional random variable with known symmetric positive definite covariance matrix $\Sigma = \mathbb{E}(\bm{\xi} \bm{\xi}^T).$  
That is, $\bm{\xi}: (\Omega, \mathcal{A}) \rightarrow (\R^n, \mathcal{B}_{\R^n})$ where $\mathcal{B}_{\R^n}$ is the Borel $\sigma$-algebra on $\R^n$. Then $\bm{\xi}$ has the joint probability density function \cite{Gauss_mult_idx_ref:master,rahman:master,xiu_thesis:master},
\begin{equation}\label{density}
\varphi: \R^n \rightarrow \R_{\ge 0} \qquad \varphi (\bm{x}; \Sigma) = \frac{1}{(2 \pi)^{n/2} |\det(\Sigma)|^{1/2} } \exp\left( - \frac{1}{2} \bm{x}^T \Sigma^{-1} \bm{x} \right).
\end{equation}

The following definition of the multivariate Hermite polynomial is a natural extension of its one-dimensional counterpart given by (\ref{1d_herm_density}). 
\begin{definition}[Multivariate Hermite polynomial]\label{mult_herm_def} Let $\bm{\xi} \sim \mathcal{N}(\bm{0}, \Sigma)$ such that $\Sigma \in \R^{n \times n}$ is symmetric positive definite (SPD) with joint density function given by Eq.~\ref{density}. Then for any $n$-dimensional multi-index $\balpha$, the \emph{multivariate Hermite polynomial} $H_{\balpha}(\bm{\xi}; \Sigma)$ indexed by $\balpha$  is a polynomial in $\bm{\xi}$ of degree $|\balpha|$ defined as 
\[ H_{\balpha}(\bm{\xi}; \Sigma) = \begin{cases} \frac{(-1)^{|\balpha|}}{\varphi(\bm{\xi}; \Sigma)} \partial^{\balpha}(\varphi(\bm{\xi}; \Sigma)) &\text{ if } |\balpha| \ge 1 \\ 1 &\text{ if } |\balpha| = 0. \end{cases} \]
We will often abbreviate a multivariate Hermite polynomial as an MHP. 
\end{definition}
Since each $H_{\balpha}$ is a polynomial of order~$\oalpha$,
\begin{align}\label{deriv_too_big_0}
\partial^{\bbeta} (H_{\balpha}) = 0 \qquad \text{whenever } \oalpha < \obeta 
\end{align}
for $\balpha, \bbeta \in \N_0^n.$ Equation (\ref{deriv_too_big_0}) will be useful in later sections.

\begin{figure}[htbp]
\centering
\includegraphics[width=0.6\textwidth]{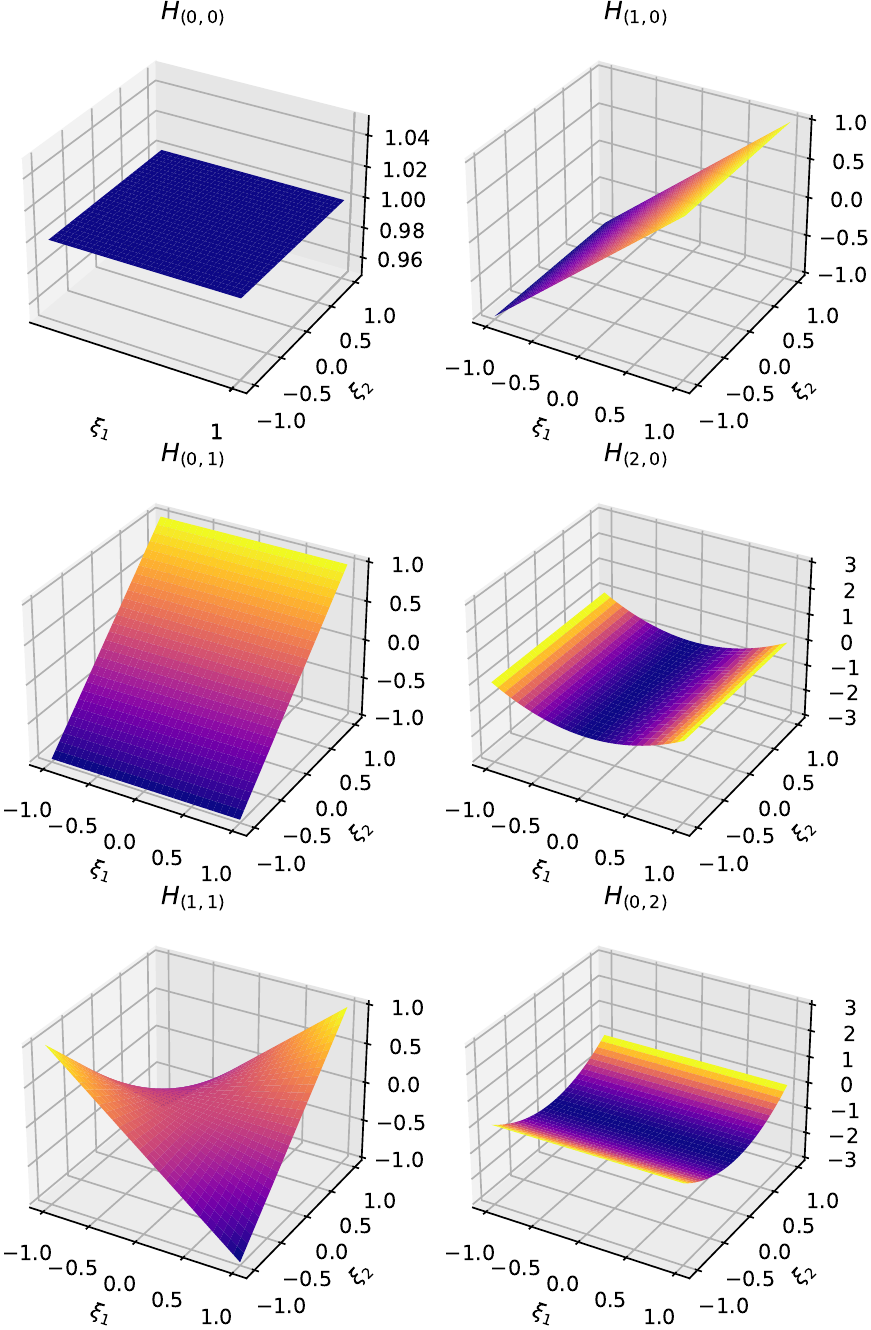}
\caption{Multivariate Hermite polynomials $H_{\balpha}(\xi_1, \xi_2)$ for $|\balpha| \le 2$, where $\mathrm{Var}(\xi_1) = \mathrm{Var}(\xi_2) = 1$ and $\mathrm{Cov}(\xi_1, \xi_2) = 1/2$.}\label{MHPsfig}
\end{figure}

Figure \ref{MHPsfig} depicts the first six MHPs graphically for a given $2 \times 2$ covariance matrix. 

The first-order moments of the MHPs are
\begin{equation}\label{firstordermoments}
\mathbb{E}(H_{\balpha}(\bxi)) = \begin{cases} 1 &\text{ if } \balpha = (0,\ldots, 0) \\ 0 &\text{ else,} \end{cases} 
\end{equation} 
as established in \cite{rahman:master}. 

For second-order moments, Dunkl and Xu \cite{dunkl_xu_2014:master} first showed that 
\begin{equation}\label{weak_orthog}
\langle H_{\balpha}, H_{\bbeta} \rangle = \mathbb{E} (H_{\bm{\alpha}}(\bm{\xi}) H_{\bm{\beta}}(\bm{\xi})) = \int H_{\balpha}(\bm{x}) H_{\bbeta}(\bm{x}) \varphi(\bm{x}) d \bm{x} = 0  \quad \text{ if } |\bm{\alpha}| \neq |\bm{\beta}|,
\end{equation}
where we recall that the double product notation $\langle H_{\balpha}, H_{\bbeta} \rangle$ was introduced by Definition \ref{inner_prod}. Rahman \cite{rahman:master} coined the condition given by (\ref{weak_orthog}) as \emph{weak orthogonality} in the context of polynomial chaos. 

By definition of weak orthogonality, $H_{\balpha}$ is orthogonal to all MHPs of different degrees but may \emph{not} be orthogonal to other MHPs of the \emph{same} $\oalpha$ degree; in fact, most multi-indices of the same degree will yield non-orthogonal MHPs, necessitating explicit double product formulas to be inroduced \cite{rahman:master, mydoubprodpaper}. 

For a given $n$ and $|\balpha| = M$, there are exactly $\mbinom{M + n -1}{n - 1}$ possible $\bbeta$ such that $|\balpha| = |\bbeta|$.\footnote{This follows from a standard \emph{stars and bars} counting argument \cite[Sect.~1.2]{Stanley2001:master}.}

When $\Sigma$ is diagonal, it is well-established \cite{Feinberg, GI_book, salazar,Smith, xiu_thesis:master} that we can decompose any $H_{\balpha}$ as the product of univariate Hermite polynomials:
\begin{align}\label{1Dprod}
H_{\balpha}(\bxi) = \prod_{i=1}^n H_{\alpha_i}(\xi_i) \quad \text{for } \bxi \sim \mathcal{N}(\bm{0}, \text{diag}(\Sigma_{11}, \ldots, \Sigma_{nn})). 
\end{align} 
Equation \ref{1Dprod} is a direct result of the PDF $\varphi$ of $\bm{\xi}$ decomposing as the product of marginal densities $\varphi_i$ for each $\xi_i$ random variable when $\bm{\xi}$ has uncorrelated components. This is what is meant by constructing multivariate orthogonal polynomials via tensoriziation of univariate orthogonal polynomials. Further, the $\{H_{\balpha}\}$ are \emph{orthogonal} in this case (not just weakly orthogonal), because
\begin{align*}
\langle H_{\balpha}, H_{\bbeta} \rangle &= \int H_{\balpha}(\bm{x}) H_{\bbeta}(\bm{x}) \varphi(\bm{x}) d \bm{x} = \int \prod_{i=1}^n H_{\alpha_i}(x_i) H_{\beta_i}(x_i) \varphi_i(x_i) d x_i \\
&= \prod_{i=1}^n \int H_{\alpha_i}(x_i) H_{\beta_i}(x_i) \varphi(x_i) dx_i = \prod_{i=1}^n \langle H_{\alpha_i}, H_{\beta_i} \rangle = \balpha! \thinspace \delta_{\balpha, \bbeta}
\end{align*}
where we applied Fubini's theorem \cite{Amir:master,Baby_Rudin} and the orthogonality of the one-dimensional $H_{\alpha_i}$ polynomials (Eq.~\ref{double_prod_1D}). When $\bm{\xi} \sim \mathcal{N}(\bm{0}, I)$, we can also observe that
\begin{align*}
H_{(0,0)} = 1, \qquad H_{(1,0)}(\xi_1, \xi_2) = \xi_1, \qquad H_{(2,0)}(\xi_1, \xi_2) = \xi_1^2 - 1, \qquad \ldots
\end{align*} 
matching the one-dimensional polynomials in Table \ref{herm_table}.

\section{Polynomial chaos with correlated Gaussians}\label{sec_PC_corr_gauss}
The MHPs serve as the basis polynomials in the multivariate polynomial chaos expansion with correlated Gaussians established by Rahman \cite{rahman:master}. In particular, we consider the following spectral series. 

\begin{definition}[Correlated multivariate polynomial chaos expansion]\label{MPC} Let $\bxi \sim \mathcal{N}(\bm{0}, \Sigma)$ and $y(\bm{\xi})$ be a real-valued output random variable defined on the same complete probability space $(\Omega, \mathcal{A}, \mathbb{P})$. If $y(\bxi) \in L^2(\Omega, \mathcal{A}, \mathbb{P})$ (that is, if $y(\bxi)$ is \emph{square-integrable}), we can express $y(\bm{\xi})$ via the \emph{correlated multivariate polynomial chaos expansion} as
\begin{equation*}
\bm{y}(\bm{\xi}) = \sum_{\balpha \in \N_0^n} y_{\balpha} H_{\balpha}(\bxi).
\end{equation*}
\end{definition}
Definition \ref{MPC} \emph{generalized} the previous notion of a multivariate PC expansion by introducing basis functionals $H_{\balpha}(\bxi)$ whose formulae incorporate the non-diagonal entries of the covariance matrix of $\bxi$.\footnote{Consequently, Rahman refers to Def.~\ref{MPC} in his original paper as the \emph{generalized PCE} \cite{rahman:master}. We instead use the name \emph{correlated PC Expansion}, or \emph{CPCE}, to avoid confusion with the \emph{generalized PCE} coined by Xiu and Karniadakis \cite{xiu_thesis:master,xiu_karniadakis:master}, which is a series that still requires independent random inputs but accounts for non-Gaussian variables.} Going forward, we will sometimes refer to the series in Definition \ref{MPC} as the \emph{correlated multivariate PCE}, the \emph{correlated PCE}, or the \emph{CPCE} for short. When referring to the multivariate PCE in \emph{independent} Gaussian parameters, we will use the phrase \emph{classical PCE}. 

When truncating to some $M$ value, we can define

\begin{equation}\label{trunc_CPCE}
y^{(M)}(\bxi) =  \sum_{\substack{\balpha \in \N_0^n \\ |\balpha| \le M}} y_{\balpha}^{(M)} H_{\balpha} (\bm{\xi}),
\end{equation}
which is a projection of $y$ onto the space $P_{M}(\bxi)$ of all polynomials  of degree up to and including $M$ in the $n$ random variables $\xi_1, \ldots, \xi_n$\cite{Smith}. 

Then $y^{(M)}(\bxi)$ converges to the true $y(\bxi)$ in mean-square, probability, and distribution~\cite{rahman:master}. The random nodes $y_{\balpha}^{(M)}$ have the superscript $(M)$ to indicate that they solve the deterministic $M$ system and thus depend on $M$ (though converge to the true $y_{\balpha}$ in Def.~\ref{MPC} as $M \rightarrow \infty$). 

The number of terms in the truncated CPCE (\ref{trunc_CPCE}) is the same as that of the classical PCE for uncorrelated Gaussians --- the only difference is that the basis functionals $H_{\balpha}(\bxi)$ now have additional (potentially nonzero) terms whose coefficients are non-diagonal entries of $\Sigma.$ Namely, the total number of terms $n_M$ in (\ref{trunc_CPCE}) is given by
\begin{align}\label{n_M}
n_M:=  \big| \{ \balpha \in \N_0^n \mid 0 \le \oalpha \le M \} \big| = \frac{(M+n)!}{M!n!} = \mbinom{M + n}{n},
\end{align}
as has been noted extensively in the literature \cite{GI_book,salazar,Smith,xiu_thesis:master,xiu_karniadakis:master}. Consequently, we can observe that correlated polynomial chaos suffers from the same \emph{curse of dimensionality} as its uncorrelated counterpart, which is to say that the number of expansion terms grows combinatorially when the dimensionality increases \cite{constantine_thesis:master,Smith,xiu_thesis:master}.

Accordingly, we return to the solving Eq.~(\ref{du_dt}), where we now assume the parameter vector is normally distributed:
\begin{align}\label{du_dt_high_dim}
 \mathcal{F}(u, \bm{d}; \bxi) &=  0, \qquad \bm{d} \in \mathcal{D} \subset \R^p, \ \bxi \sim \mathcal{N}(\bm{0}, \Sigma).
 \end{align}

If $u(\bm{d}; \bxi) \in L^2(\Omega, \mathcal{A}, \mathbb{P})$ for every $\bm{d} \in \mathcal{D}$, we can approximate $u(\bm{d}; \bxi)$ using the truncated CPCE:
\begin{equation}\label{high_dim_PCE_old}
u(\bm{d}; \bm{\xi}) = \sum_{\balpha \in \N_0^n} \thinspace u_{\balpha}(\bm{d}) H_{\balpha}(\bm{\xi})
\end{equation}
where the basis functionals $\Psi_{\balpha}$ from (\ref{du_dt}) are now chosen to be multivariate Hermite polynomials (Def.~\ref{mult_herm_def}).



\section{Example application: logistic growth modeling}\label{example_application}

Consider the following stochastic variant of Hutchinson's equation \cite{hutch}, which models population growth with a time lag,
\begin{align}\label{logex}
 \partial_t X(t;\bxi) &= r(\bxi) X(t; \bxi) \big[ 1 - K X(t - \tau; \bxi) \big]  \\
 X(t;\bxi) &= 0 \text{ for } t \le 0.
 \end{align}
In \Cref{logex}, 
\begin{itemize}
\item the solution $X(t; \bxi)$ is the population size at time $t$, 
\item $r(\bxi)$ is the stochastic growth rate (whose values are real and possibly negative), 
\item $\tau > 0$ is a constant discrete delay, and
\item $K > 0$ is the constant carrying capacity. 
\end{itemize}
We assume $X(t;\bxi), r(\bxi) \in L^2(\Omega, \mathcal{A}, \mathbb{P})$, and $\bxi \sim \mathcal{N}(\bm{0}_{n}, \Sigma)$ for known SPD covariance matrix $\Sigma \in \R^{n \times n}$. The truncated CPCE for a given truncation parameter $M \ge 0$ gives

\begin{align*}
X(t;\bxi) &\approx \sum_{|\bm{i}| \le M} x_{\bm{i}}(t) H_{\bm{i}}(\bxi)  \\
r(\bxi) &\approx \sum_{|\bm{i}| \le M} r_i H_{\bm{i}}(\bxi) 
\end{align*} 
By substituting these approximations into \Cref{logex}, we have
\begin{align*}
\sum_{|\bm{i}| \le M}  \partial_t x_{\bm{i}}(t) H_{\bm{i}}(\bxi) &= \sum_{|\bm{i}|, |\bm{j}| \le M} r_{\bm{i}} x_{\bm{j}}(t) H_{\bm{i}}(\bxi) H_{\bm{j}}(\bxi) \\
&- K   \sum_{|\bm{i}|, |\bm{j}|, |\bm{k}| \le M} r_{\bm{i}} x_{\bm{j}}(t) x_{\bm{k}}(t - \tau) H_{\bm{i}}(\bxi) H_{\bm{j}}(\bxi)  H_{\bm{k}}(\bxi).
\end{align*}
We multiply through by $H_{\bm{\ell}}$ for fixed but arbitrary $\bm{\ell} \in \N_0^n, 0 \le |\bm{\ell}| \le M$ and integrate (i.e., the Galerkin projection step). Then we have the system of deterministic PDEs

\begin{equation}\label{example_sys}
 \begin{cases} \sum_{|\bm{i}| = |\bm{\ell}|}  \partial_t x_{\bm{i}}(t) \langle H_{\bm{i}}, H_{\bm{\ell}} \rangle &= \sum_{|\bm{i}|, |\bm{j}| \le M} r_{\bm{i}} x_{\bm{j}}(t) \langle H_{\bm{i}}, H_{\bm{j}}, H_{\bm{\ell}} \rangle  \\
&- K   \sum_{|\bm{i}|, |\bm{j}|, |\bm{k}| \le M} r_{\bm{i}} x_{\bm{j}}(t) x_{\bm{k}}(t - \tau) \langle H_{\bm{i}}, H_{\bm{j}}, H_{\bm{k}}, H_{\bm{\ell}} \rangle
 \end{cases} 
 \end{equation}
 for $0 \le |\bm{\ell}| \le M.$ Note that for the summation on the left-hand side of (\ref{example_sys}) we have used weak orthogonality to eliminate the $\bm{i}$ such that $|\bm{i}| \neq |\bm{\ell}|$. 
 
We can quickly see the \emph{large number of double, triple, and quadruple product coefficients} that need to be calculated in order to even express (and then solve) the above PDE system. For concrete examples, consider the following.

\begin{itemize}
\item Suppose $M = 2$ and $n = 2$. Then the system of PDEs in (\ref{example_sys}) has ${M + n \choose n}$ = 6 equations total (by Equation \ref{n_M}), one for each value of $\bm{\ell}$. Within each equation, the summation over the $\bm{i}$ such that $|\bm{i}| = |\bm{\ell}|$ will have $\mbinom{M + n -1}{n - 1} = 3$ terms. The summation over $ 0 \le |\bm{i}|, |\bm{j}| \le M$ gives at most 12 nonzero terms, and the summation over $0 \le |\bm{i}|, |\bm{j}|, |\bm{\ell}| \le M$ gives at most 18 nonzero terms. Therefore, each equation has up to $3 + 12 + 18 = 33$ $N$th products. Since there are 6 equations in total, this amounts to at most $33 \times 6 = 189$ coefficients to compute.

\item Suppose $M = 2$ and $n = 3$. Then the system of PDEs in (\ref{example_sys}) has ${M + n \choose n}$ = 10 equations total. Within each equation, the summation over the $\bm{i}$ such that $|\bm{i}| = |\bm{\ell}|$ will have $\mbinom{M + n -1}{n - 1} = 6$ terms. The summation over $ 0 \le |\bm{i}|, |\bm{j}| \le M$ gives at most 20 terms, and the summation over $0 \le |\bm{i}|, |\bm{j}|, |\bm{\ell}| \le M$ gives at most 30 nonzero terms.
\item Therefore, each equation has up to $6 + 20 + 30 = 56$ $N$th products, and since there are 10 equations in total, this amounts to an upper bound of 560 coefficients to compute. 
\end{itemize}

Consequently, we can see that estimating each \emph{individual} $N$th product coefficient via Monte Carlo sampling can quickly become cost prohibitive, especially if each \emph{individal} evaluation takes on the order of (say) minutes. Therefore, having a closed-form equation for the $N$th products has the potential to  nontrivially reduce total computation time resulting from the curse of dimensionality in SGS.

\section{Formulae for multivariate Hermite polynomials}\label{multivar_formulae}
This section overviews decompositions of MHPs and their derivatives previously established \cite{mytripprodpaper}. Namely, \cref{Nthprodformula} combines them with the integration by parts identity of \cref{IBP} (proved in \cref{existence_section}) to obtain the $N$th product recursion. For ease of notation, we sometimes omit the argument $\bxi$ of an MHP $H_{\balpha}$, and all expectations involving MHPs are taken with respect to the law of $\bxi \sim \mathcal{N}(\bm{0}, \Sigma)$.

In \cite{mydoubprodpaper}, we established the following result.



We begin by rewriting a MHP indexed by a sum of indices as a finite sum of lower-order Hermite polynomials multiplied by certain partial derivatives. 
\begin{lemma}[Sum of indices formula \cite{mytripprodpaper}]\label{sum_of_indices}
Let $\balpha, \bbeta \in \N_0^n.$ Then
\[ H_{\balpha + \bbeta} = \sum_{\bm{q} \le \balpha} (-1)^{|\bm{q}|} \mbinom{\balpha}{\bm{q}} H_{\balpha - \bm{q}} \thinspace \partial^{\bm{q}} (H_{\bbeta}). \] 
\end{lemma}


\begin{theorem}[Higher order derivative of MHP \cite{mytripprodpaper}]\label{deriv_formula}Let $\A, \bbeta \in \N_0^n$ such that $|\A| \le \obeta.$ Then 
\[ \partial^{\A} (H_{\bbeta}) = \sum_{\substack{\B \le \bbeta \\ |\B| = |\A|}} \mbinom{\bbeta}{\B}  \langle H_{\A}, H_{\B} \rangle \thinspace H_{\bbeta - \B}. \] 
\end{theorem}



Theorem \ref{deriv_formula} also helps tie up a technical and \emph{subtle} loose end. Before, Theorem \ref{deriv_formula} established that $\oalpha = \obeta$ implies $\mathbb{E}(\partial^{\balpha} H_{\bbeta}) = \doub.$ It turns out that $\mathbb{E}(\partial^{\balpha} H_{\bbeta}) = \doub$ in \emph{general}, i.e., even when  $\oalpha \neq \obeta.$ By \cite{rahman:master}, we of course know that $\doub = 0$ when $\oalpha \neq \obeta$. By proving  $\mathbb{E}(\partial^{\balpha} H_{\bbeta}) = 0$ when $\oalpha \neq \obeta$, we will have shown that $\mathbb{E}(\partial^{\balpha} H_{\bbeta}) = \doub$ for any $\balpha, \bbeta \in \N_0^n.$ 

Clearly $\partial^{\balpha} H_{\bbeta} = 0$ when $\oalpha > \obeta$, since $H_{\bbeta}$ is a polynomial of order $\obeta$, i.e., (\ref{deriv_too_big_0}); thus, $\mathbb{E}(\partial^{\balpha} H_{\bbeta})  = 0$ in this case. When $\oalpha < \obeta,$ \Cref{deriv_formula} gives us terms each with $\mathbb{E}(H_{\bbeta - \B})$ such that $\bbeta - \B \neq \bm{0}$ (since $|\B| = \oalpha < \obeta$); hence, we have $\partial^{\balpha} (H_{\bbeta}) = 0$ in this case as well. Consequently, the general derivative formula (Theorem \ref{deriv_formula}) establishes that $\mathbb{E}(\partial^{\balpha} H_{\bbeta}) = \doub$ for any multi-indices $\balpha, \bbeta \in \N_0^n.$

\section{The $N$th product formula}\label{Nthprodformula}

\begin{theorem}[$N$th order product]\label{N_prod} Let $\balpha, \bbeta, \bgamma^{(1)}, \ldots, \bgamma^{(N-2)} \in \N_0^n$ such that $3 \le N <~\infty$ and $\oalpha \le \obeta \le |\bgamma^{(1)}| \le \cdots \le |\bgamma^{(N-2)}|$. Then
\begin{align*}
&\langle H_{\balpha}, H_{\bbeta}, H_{\bgamma^{(1)}}, \ldots, H_{\bgamma^{(N-2)}} \rangle \\
&= \sum_{\A \le \balpha} \sum_{\substack{\B \le \bbeta \\ |\B| = |\A|}} \mbinom{\balpha}{\A} \mbinom{\bbeta}{\B} \langle H_{\A}, H_{\B} \rangle \langle H_{\balpha + \bbeta - \A - \B}, H_{\bgamma^{(1)}}, \ldots, H_{\bgamma^{(N-2)}} \rangle.
\end{align*}
\end{theorem}
\begin{proof}
Observe that
\begin{flalign*} &\hphantom{=} \langle H_{\balpha}, H_{\bbeta}, H_{\bgamma^{(1)}}, \ldots, H_{\bgamma^{(N-2)}} \rangle \\
&= \mathbb{E}_{\bxi}[\partial^{\balpha} (H_{\bbeta}(\bxi) H_{\bgamma^{(1)}}(\bxi) \cdots H_{\bgamma^{(N-2)}}(\bxi)) &\text{[Corollary \ref{general_deriv_N_product}]} \\
&= \int_{\R^n} \partial^{\balpha}(H_{\bbeta}(\x) H_{\bgamma^{(1)}}(\x) \cdots H_{\bgamma^{(N-2)}}(\x) )\varphi(\x) \thinspace d \x \\
&= \sum_{\A \le \balpha} \mbinom{\balpha}{\A} \int \partial^{\A} (H_{\bbeta}(\bm{x})) \partial^{\balpha - \A}(\textstyle \prod_{i=1}^{N-2} H_{\bgamma^{(i)}}(\x) ) \varphi(\bm{x}) \thinspace d \bm{x}   \hspace{50pt} & \text{[Eq.~\ref{leibniz-rule}]}   \\
&= \sum_{\A \le \balpha} \sum_{\substack{\B \le \bbeta \\ |\B| = |\A|}} \mbinom{\balpha}{\A} \mbinom{\bbeta}{\B} \langle H_{\A}, H_{\B} \rangle \hspace{-4pt}  \int \hspace{-4pt} H_{\bbeta - \B}(\bm{x}) \thinspace \partial^{\balpha - \A}\left(\textstyle \prod_{i=1}^{N-2} H_{\bgamma^{(i)}}(\x) \right) \varphi(\bm{x}) \thinspace d \bm{x} \hspace{-55pt} 
\end{flalign*}
where the last equality follows from \Cref{deriv_formula}. Continuing the calculation,

\begin{flalign*}
&\langle H_{\balpha}, H_{\bbeta}, H_{\bgamma^{(1)}}, \ldots, H_{\bgamma^{(N-2)}} \rangle \\
&= \sum_{\A \le \balpha} \sum_{\substack{\B \le \bbeta \\ |\B| = |\A|}} \mbinom{\balpha}{\A} \mbinom{\bbeta}{\B} \langle H_{\A}, H_{\B} \rangle \thinspace \mathbb{E}_{\bxi} [ H_{\bbeta - \B}(\bxi) \partial^{\balpha - \A}(\textstyle \prod_{i=1}^{N-2} H_{\bgamma^{(i)}}(\bxi) ) ] \\
&= \sum_{\A \le \balpha} \sum_{\substack{\B \le \bbeta \\ |\B| = |\A|}} \mbinom{\balpha}{\A} \mbinom{\bbeta}{\B} \langle H_{\A}, H_{\B} \rangle \thinspace \mathbb{E}_{\bxi} \left[ \partial^{\balpha + \bbeta - \A - \B} (\textstyle \prod_{i=1}^{N-2} H_{\bgamma^{(i)}}(\bxi) ) \right]  & \text{[Thm.~\ref{IBP}]} \\
&= \sum_{\A \le \balpha} \sum_{\substack{\B \le \bbeta \\ |\B| = |\A|}} \mbinom{\balpha}{\A} \mbinom{\bbeta}{\B} \langle H_{\A}, H_{\B} \rangle \langle H_{\balpha + \bbeta - \A - \B}, H_{\bgamma^{(1)}}, \ldots, H_{\bgamma^{(N-2)}} \rangle. & \text{[Corollary \ref{general_deriv_N_product}]} 
\end{flalign*}
Theorem \ref{IBP} could be applied in the second-to-last equality above, since Lemma \ref{slowly-growing-lemma} ensures that $\partial^{\balpha - \A} (\textstyle \prod_{i=1}^{N-2} H_{\bgamma^{(i)}})~\in~\mathcal{S}(\R^n)$. 
\end{proof}

\begin{corollary}[Zeros of $N$th order product]\label{N_prod_zeros}Let $\balpha, \bbeta, \bgamma^{(1)}, \ldots, \bgamma^{(N-2)} \in \N_0^n$ such that $3 \le N <~\infty$ and $\oalpha \le \obeta \le |\bgamma^{(1)}| \le \cdots \le |\bgamma^{(N-2)}|$. If 

\begin{align*}
 \oalpha + \obeta + |\bgamma^{(1)}| + \cdots + |\bgamma^{(N-3)}|  < |\bgamma^{(N-2)}| \tag*{for $N > 3$,}
 \end{align*}
or
\begin{align*}
\oalpha + \obeta < |\bgamma^{(1)}|  \tag*{for $N = 3$,}
\end{align*}
then $\langle H_{\balpha}, H_{\bbeta}, H_{\bgamma^{(1)}}, \ldots, H_{\bgamma^{(N-2)}} \rangle = 0.$
\end{corollary} 
\begin{proof}We can see this result in one of two ways: by applying Corollary \ref{general_deriv_N_product} directly, or using the formula from Theorem \ref{N_prod} above and then applying Corollary \ref{general_deriv_N_product} to each term in the summation. For the former, Corollary \ref{general_deriv_N_product} gives\footnote{Note that if $N = 3$, we interpret the product in (\ref{eq_in_n_prod_0_arg}) as empty (and thus equal to 1), meaning we have $\langle H_{\balpha}, H_{\bbeta}, H_{\bgamma^{(1)}} \rangle =  \mathbb{E}_{\bxi}\left(\partial^{\bgamma^{(1)}}\left(H_{\balpha}(\bxi) H_{\bbeta}(\bxi) \right) \right).$ }
\begin{align}\label{eq_in_n_prod_0_arg}
\langle H_{\balpha}, H_{\bbeta}, H_{\bgamma^{(1)}}, \ldots, H_{\bgamma^{(N-2)}} \rangle = \mathbb{E}_{\bxi}\left(\partial^{\bgamma^{(N-2)}} \left( H_{\balpha}(\bxi) H_{\bbeta}(\bxi) \textstyle\prod_{i=1}^{N-3} H_{\bgamma^{(i)}}(\bxi)\right) \right).
\end{align}
However, note that $H_{\balpha} H_{\bbeta} \prod_{i=1}^{N-3} H_{\bgamma^{(i)}}$ is a polynomial of order $\oalpha + \obeta + |\bgamma^{(1)}| + \cdots + |\bgamma^{(N-3)}|$. So differentiating by the multi-index $\bgamma^{(N-2)}$ with $|\bgamma^{(N-2)}| > \oalpha + \obeta + |\bgamma^{(1)}| + \cdots + |\bgamma^{(N-3)}|$ gives the zero polynomial in $\bxi$, whose expected value is then zero as well. 

Alternatively, we could use the formula for the $N$th order product from Theorem \ref{N_prod}. In this case, we apply the same argument just given (via Corollary \ref{general_deriv_N_product}) to show that every term contains a factor
\[ \langle H_{\balpha + \bbeta - \A - \B}, H_{\bgamma^{(1)}}, \ldots, H_{\bgamma^{(N-2)}} \rangle = 0. \]
\end{proof}



\subsection{Examples of quadruple product computations}\label{sect_trip_prod_calcs}

\begin{example}\label{ex_3} For $n = 3,$ let $\bm{\alpha} = (1,0,2), \bm{\beta} = (3,1,0),$ and $\bm{\gamma} = (1,2,2).$ Then $|\bm{\alpha}| = 3,$ $|\bm{\beta}| = 4,$ and $|\bm{\gamma}| = 5,$ so we cannot use conditions on size or parity to conclude the triple product is zero. In this case, $k := \frac{|\bm{\alpha}| + |\bm{\beta}| - |\bm{\gamma}|}{2} = 1.$

We need to find the $\bm{A} \le \bm{\alpha}$ and $\B \le \bbeta$ such that $|\bm{A}| = |\B| = 1$, so the options are $\A = \bm{e_1}, \bm{e_3}$ and $\B = \bm{e_1}, \bm{e_2}.$ The corresponding terms are computed in Table  \ref{tab_3}.

\begin{table}[htbp]
\caption{Terms of the triple product formula for $\bm{\alpha} = (1,0,2)$, $\bm{\beta} = (3,1,0)$, and $\bm{\gamma} = (1,2,2)$.}\label{tab_3}
\centering
\renewcommand{\arraystretch}{1.2}
\begin{tabular}{ccccccc}
\toprule
$\bm{A}$ & $\bm{B}$ & $\mbinom{\balpha}{\A}$ & $\mbinom{\bbeta}{\B}$ & $\langle H_{\A}, H_{\B} \rangle$ & $H_{\bm{\alpha + \beta - A - B}}$ & Term \\
\midrule
$\bm{e_1}$ & $\bm{e_1}$ & 1 & 3 & $\Sigma_{11}^{-1}$ & $H_{(2,1,2)}$ & $3 \Sigma_{11}^{-1} \langle H_{(2,1,2)}, H_{\bm{\gamma}} \rangle$ \\
$\bm{e_1}$ & $\bm{e_2}$ & 1 & 1 & $\Sigma_{12}^{-1}$ & $H_{(3,0,2)}$ & $\Sigma_{12}^{-1} \langle H_{(3,0,2)}, H_{\bm{\gamma}} \rangle$ \\
$\bm{e_3}$ & $\bm{e_1}$ & 2 & 3 & $\Sigma_{13}^{-1}$ & $H_{(3,1,1)}$ & $6 \Sigma_{13}^{-1} \langle H_{(3,1,1)}, H_{\bm{\gamma}} \rangle$ \\
$\bm{e_3}$ & $\bm{e_2}$ & 2 & 1 & $\Sigma_{23}^{-1}$ & $H_{(4,0,1)}$ & $2 \Sigma_{23}^{-1} \langle H_{(4,0,1)}, H_{\bm{\gamma}} \rangle$ \\
\bottomrule
\end{tabular}
\end{table}

Therefore,
\begin{align*}
\langle H_{(1,0,2)}, H_{(3,1,0)}, H_{(1,2,2)} \rangle &=  \sum_{\substack{\A \le (1,0,2) \\ \B \le (3,1,0) \\ |\A| = |\B| = 1}} \mbinom{\balpha}{\A} \mbinom{\bbeta}{\B} \langle H_{\A}, H_{\B} \rangle \langle H_{\bm{\alpha + \beta - \A - \B}}, H_{\bgamma} \rangle \\
&= 3 \Sigma_{11}^{-1} \langle H_{(2,1,2)}, H_{(1,2,2)} \rangle + \Sigma_{12}^{-1} \langle H_{(3,0,2)}, H_{(1,2,2)} \rangle  \\
&+ \thinspace 6 \Sigma_{13}^{-1} \langle H_{(3,1,1)}, H_{(1,2,2)} \rangle + 2 \Sigma_{23}^{-1} \langle H_{(4,0,1)}, H_{(1,2,2)} \rangle. 
\end{align*}
\end{example}


\begin{example}\label{ex_gamma_too_big} For $n = 4,$ let $\bm{\alpha} = (1,0,0,1),$ $\bm{\beta} = (0,2,1,0)$, $\bgamma^{(1)} = (0,1,1,1)$, and $\bgamma^{(2)} = (1,3,2,3).$ Then $\oalpha = 2$, $\obeta = 3$, $|\bgamma^{(1)}| = 3$, and $|\bgamma^{(2)}|= 9$, so $\oalpha + \obeta + |\bgamma^{(1)}| < |\bgamma^{(2)}|.$ We conclude Corollary \ref{N_prod_zeros} that $\langle H_{(1,0,0,1)}, H_{(0,2,1,0)}, H_{(0,1,1,1)}, H_{(1,3,2,3)} \rangle = 0.$ 
\end{example}

\begin{example}\label{quad_prod_calc} For $n = 2$, let $\balpha = (1,0)$, $\bbeta = (1,1)$, $\bgamma^{(1)} = (1,1)$, and $\bgamma^{(2)} = (3,2).$  Then $\oalpha + \obeta + |\bgamma^{(1)}| =  |\bgamma^{(2)}|$, so we cannot  immediately apply Corollary \ref{N_prod_zeros} to conclude that the quadruple product is zero. By Theorem \ref{N_prod},

\begin{align*}
\langle H_{\balpha}, H_{\bbeta}, H_{\bgamma^{(1)}}, H_{\bgamma^{(2)}} \rangle 
&= \sum_{\A \le (1,0)} \sum_{\substack{\B \le (1,1) \\ |\B| = |\A|}} \mbinom{\balpha}{\A} \mbinom{\bbeta}{\B} \langle H_{\A}, H_{\B} \rangle \langle H_{\balpha + \bbeta - \A - \B}, H_{\bgamma^{(1)}}, H_{\bgamma^{(2)}} \rangle.
\end{align*}
Table \ref{tab_4} calculates each of the terms in this summation. 
 
\begin{table}[htbp]
\caption{Terms of the quadruple ($N = 4$) product formula of \cref{N_prod} for $\bm{\alpha} = (1,0)$, $\bm{\beta} = \bm{\gamma}^{(1)} = (1,1)$, and $\bgamma^{(2)} = (3,2)$.}\label{tab_4}
\centering
\renewcommand{\arraystretch}{1.2}
\begin{tabular}{ccccccc}
\toprule
$\bm{A}$ & $\bm{B}$ & $\mbinom{\balpha}{\A}$ & $\mbinom{\bbeta}{\B}$ & $\langle H_{\A}, H_{\B} \rangle$ & $H_{\bm{\alpha + \beta - A - B}}$ & Term \\
\midrule
$\bm{0}$ & $\bm{0}$ & 1 & 1 & 1 & $H_{(2,1)}$ & $\langle H_{(1,1)}, H_{(2,1)}, H_{(3,2)} \rangle$ \\
$\bm{e_1}$ & $\bm{e_1}$ & 1 & 1 & $\Sigma_{11}^{-1}$ & $H_{(0,1)}$ & $\Sigma_{11}^{-1} \langle H_{(0,1)}, H_{(1,1)}, H_{(3,2)} \rangle$ \\
$\bm{e_1}$ & $\bm{e_2}$ & 1 & 1 & $\Sigma_{12}^{-1}$ & $H_{(1,0)}$ & $\Sigma_{12}^{-1} \langle H_{(1,0)}, H_{(1,1)}, H_{(3,2)} \rangle$ \\
\bottomrule
\end{tabular}
\end{table}
\noindent Thus, 
\begin{align*}
\langle H_{(1,0)}, H_{(1,1)}, H_{(1,1)}, H_{(3,2)} \rangle &= \langle H_{(1,1)}, H_{(2,1)}, H_{(3,2)} \rangle + \Sigma_{11}^{-1} \langle H_{(0,1)}, H_{(1,1)}, H_{(3,2)} \rangle \\
&+ \Sigma_{12}^{-1} \langle H_{(1,0)}, H_{(1,1)}, H_{(3,2)} \rangle.
\end{align*} 
However, here we can notice that $\langle H_{(0,1)}, H_{(1,1)}, H_{(3,2)} \rangle =  \langle H_{(1,0)}, H_{(1,1)}, H_{(3,2)} \rangle = 0,$ since $|(3,2)| = 5$ is greater than $|(0,1)| + |(1,1)|$ and $|(1,0)| + |(1,1)|$. Therefore,
\begin{align*}
\langle H_{(1,0)}, H_{(1,1)}, H_{(1,1)}, H_{(3,2)} \rangle &= \langle H_{(1,1)}, H_{(2,1)}, H_{(3,2)} \rangle \\
&= \langle H_{(3,2)}, H_{(3,2)} \rangle.
\end{align*}
The triple product equals the double product $\langle H_{(3,2)}, H_{(3,2)} \rangle$ above by \cite{mytripprodpaper}, since $k = \frac{|(1,1)| + |(2,1)| - |(3,2)| }{2}$  $= 0.$
\end{example}

In the next section, we computationally verify moment calculations like the ones considered in Examples \ref{ex_3}, \ref{ex_gamma_too_big}, and \ref{quad_prod_calc}. In particular, we compare the moments predicted by the formulae introduced in this paper to estimates from simple Monte Carlo sampling for a variety of test cases. 

\section{Comparison with Monte Carlo estimates}\label{simple_MC_conv}
The test SPD covariance matrices are
\begin{align}\label{cov_test_matrices}
 \Sigma^{(1)} &= \begin{pmatrix} \hphantom{-}1  & - 1/2 \\ - 1/2 & \vphantom{\frac{1}{2}}\phantom{-}1 \end{pmatrix} \qquad 
 \Sigma^{(2)} =  \begin{pmatrix} 1 & \frac{1}{5}  & \frac{2}{5} \\ \frac{1}{5} & 1 & \frac{4}{5} \\ \frac{2}{5} & \frac{4}{5} & 1 \end{pmatrix}.
\end{align}
We simulate samples drawn for $\bm{\xi} \sim \mathcal{N}(\bm{0}, \Sigma)$ via \texttt{random.multivariate\_normal()} in the \texttt{numpy} package of \texttt{Python3}.  Given these realizations $\bxi^{(i)}$ for $i = 1, \ldots, \EuScript{N}$ and a given integer $\EuScript{N}$, we can estimate $ \langle H_{\balpha}, H_{\bbeta}, H_{\bgamma} \rangle = \mathbb{E}_{\bxi} (H_{\balpha}(\bxi) H_{\bbeta}(\bxi) H_{\bgamma}(\bxi))$ by the sample mean,

\begin{align}\label{S_N}
S_{\EuScript{N}} := \frac{1}{\EuScript{N}} \sum_{i=1}^{\EuScript{N}} H_{\balpha}(\bxi^{(i)})H_{\bbeta}(\bxi^{(i)}) H_{\bgamma}(\bxi^{(i)})
\end{align}
with convergence via the Central Limit Theorem (CLT) at a rate of $O(\frac{1}{\sqrt{\EuScript{N}}})$ \cite{constantine_thesis:master,Amir:master,mcbook:master, Smith}. By symmetry of the $N$th product, we choose our index labels such that $\oalpha \le \obeta \le \ogamma.$ 

For multi-indices in $\N_0^2$, we calculate he sample means for estimating the quadruple products $\langle H_{\balpha}, H_{\bbeta}, H_{\bgamma}, H_{\bm{\delta}} \rangle$ with the same range of sample $\EuScript{N}$ counts for all 70 combinations of $\balpha, \bbeta, \bgamma, \bm{\delta}$ with orders $0 < \oalpha, \obeta, \ogamma, \odelta \le 2$ (where $H_{(0,0)} = 1$ is excluded, since these reduce to single, double, or triple products). The corresponding data for all $3 \times 1$ multi-indices is available as well~\cite{mydatataable}. 
\begin{figure}[htbp]
\centering
\includegraphics[width=0.49\textwidth]{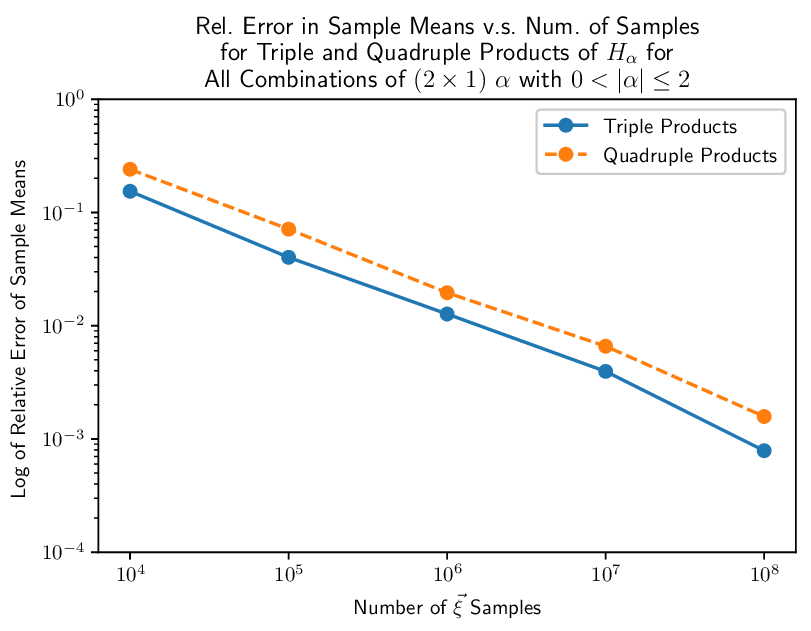}\hfill
\includegraphics[width=0.49\textwidth]{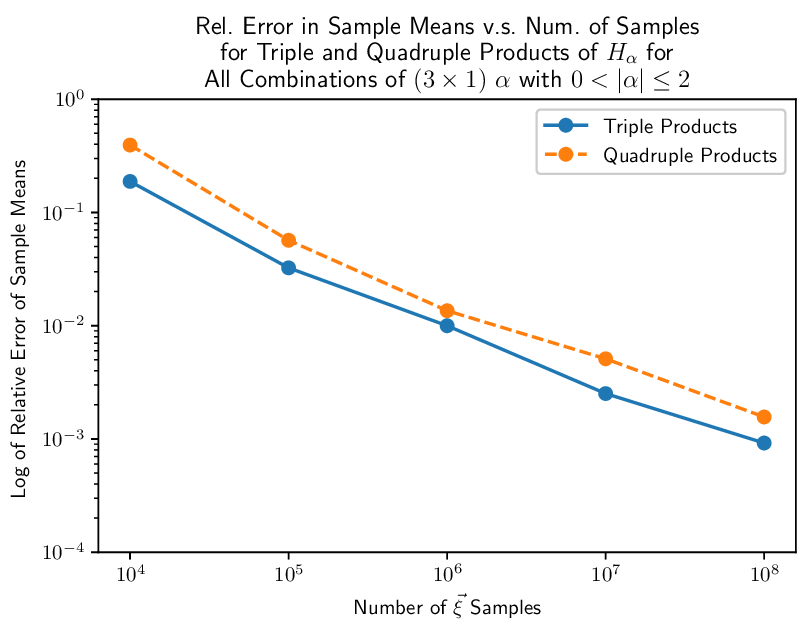}
\caption{Log of relative error for the triple and quadruple products estimated by simple Monte Carlo versus the formula from Thm.~\ref{N_prod} for $\bxi = (\xi_1, \xi_2)^T \sim \mathcal{N}(\bm{0}, \Sigma^{(1)})$ (LHS) and $\bxi = (\xi_1, \xi_2, \xi_3)^T \sim \mathcal{N}(\bm{0}, \Sigma^{(2)})$ (RHS), where the covariances $\Sigma^{(1)}$  and $\Sigma^{(2)}$ are given in (\ref{cov_test_matrices}). For the triple products, $\bm{x}_{\text{true}}$ is the vector of $\trip$ containing the 35 combinations of $\balpha, \bbeta, \bgamma \in \N_0^2$ (on the LHS) and the 165 combinations of $\balpha, \bbeta, \bgamma \in \N_0^3$ (on the RHS) such that $0 < \oalpha \le \obeta \le \ogamma \le 2$. For the quadruple products, $\bm{x}_{\text{true}}$ is the vector of $\langle H_{\balpha}, H_{\bbeta},  H_{\bgamma}, H_{\bdelta} \rangle$ containing the 70 combinations of $\balpha, \bbeta, \bgamma, \bdelta \in \N_0^2$ (on the LHS) and the 475 combinations of $\balpha, \bbeta, \bgamma, \bdelta \in \N_0^3$ (on the RHS) such that $0 < \oalpha \le \obeta \le \ogamma \le |\bm{\delta}| \le 2$. In both cases, $\bm{x}_{\text{est.}}^{(\EuScript{N})}$ is the vector 
containing the sample means $S_{\EuScript{N}}$ that estimate the corresponding products for each fixed $\EuScript{N}$ number of samples ($\EuScript{N} = 10^4, \ldots, 10^8$). Note that the $H_{\bm{0}} = 1$ polynomials are excluded, since these triple and quadruple products reduce to lower-order products.}\label{fig_relerr_n3_Nprods}
\end{figure}


Figure \ref{fig_relerr_n3_Nprods} summarizes the convergence of simple Monte Carlo to the predicted triple and quadruple product values. For each $\EuScript{N}$ sample count ($\EuScript{N} = 10^4, \ldots, 10^8$), the relative error is calculated as $\frac{\| \bm{x}_{\text{true}} - \bm{x}_{\text{est.}}^{(\EuScript{N})} \|_2}{\| \bm{x}_{\text{true}} \|_2}$. We can observe that the MC estimates conform with the predicted product values, though a large number of realizations of $\bxi$ is required --- namely, $O(10^8)$ samples for relative error to fall to around $O(10^{-3})$ for both test correlation matrices. 

\section{Conclusions}
\label{sec:conclusions}


This paper provides several theoretical results for correlated multivariate Hermite polynomials (MHPs), some of which are as follows.
\begin{itemize}
\item We prove a formula for the expected value of an MHP times a so-called \emph{slowly growing} function, both with correlated Gaussian inputs, which is used to establish the existence of the $N$th product  (\Cref{IBP,moments-exist,general_deriv_N_product}).
\item We establish an expression for the $N$th product of MHPs in correlated Gaussian random variables  (\Cref{N_prod}). 
\item We determine sufficient conditions for the $N$th product  to equal zero (Corollary \ref{N_prod_zeros}). 

\end{itemize}

In \cref{existence_section}, we  establish that the expected value with respect to the $d\bm{\xi}$ measure exists for any slowly growing statistic of the Gaussian $\bxi$ random variable (Theorem \ref{schwartz_gives_moments_exist}). This further highlights the importance of knowing definitively when $N$th products are zero. When the moments are approximated numerically, it can be difficult to distinguish between a genuine small physical effect in the $M$-system (that is, a PDE term whose coefficient is near zero) and numerical noise from the estimation process; this is especially important if such a perturbation nontrivially affects the system solution.  Therefore, it would be valuable to know a priori which products are 0 so that their corresponding terms in the $M$-system could be eliminated altogether with confidence.

Finally, the delayed logistic growth model of \cref{example_application} illustrates how quickly the number of required $N$th products grows. The numerical experiments of \cref{simple_MC_conv} indicate that MC estimates of these products approach the values given by \cref{N_prod} only after roughly $10^8$ samples each. Evaluating the formula instead removes this sampling cost from the assembly of the stochastic Galerkin system for nonlinear PDEs with correlated Gaussian parameters.


\appendix

\section{Slowly growing and Schwartz functions}\label{appendix_defs}
The following definitions of \emph{slowly growing} and \emph{Schwartz functions} \cite{slowly-growing-plus-prod-schwartz}  will assist with formalizing the types of functions $f$ for which $\mathbb{E}_{\bm{\xi}}(f(\bm{\xi}))$ is guaranteed to exist when $\bxi \sim \mathcal{N}(\bm{0}, \Sigma).$ 

\begin{definition}[Slowly growing]\label{slowly_increasing}
We say that $f \in C^\infty(\R^n, \R)$ is \emph{slowly growing} if $f$ and all of its derivatives have at most polynomial growth; that is, for every $\balpha \in \N_0^n,$ there exists a constant $C_{\balpha}$ and an integer $N_{\balpha} \in \N_0$ such that 
\[ | \partial^{\balpha} f (\bm x) | \le C_{\balpha} (1 + |\bm x|^2)^{N_{\balpha}} \qquad \qquad \text{ for all } \bm{x} \in \R^n. \] 
We will let $\mathcal{P}(\R^n)$ be the set of all slowly growing, real-valued functions on $\R^n$.
\end{definition}
 \begin{definition}[Schwartz functions]\label{Schwartz}The set $\mathcal{S}(\R^n)$, called the set of \emph{Schwartz functions}, consists of functions $\phi \in C^\infty(\R^n, \R)$ such that for all multi-indices $\balpha, \bbeta \in \N_0^n,$ \\ $\sup_{\bm{x} \in \R^n} |\bm{x}^{\balpha} \partial^{\bbeta}(\phi)| < \infty.$
\end{definition}
The relationships given in Lemma \ref{slowly-growing-lemma} below are well-established \cite{MSE:prod-schwartz-poly, slowly-growing-plus-prod-schwartz}, but we present them for use in later theorems. 
\begin{lemma}\label{slowly-growing-lemma}Let $\mathcal{P}(\R^n)$ be the set of slowly growing, real-valued functions on $\R^n$. Then:
\begin{enumerate}
\item (Closure under point-wise function multiplication) $f, g \in \mathcal{P}(\R^n) \Rightarrow f g \in \mathcal{P}(\R^n)$;
\item (Closure under differentiation) $f \in \mathcal{P}(\R^n) \Rightarrow \partial^{\balpha} f \in \mathcal{P}(\R^n) \text{ for all } \balpha \in \N_0^n$;
\item (Slowly growing times Schwartz is Schwartz) $ f \in\mathcal{P}(\R^n), \phi \in \mathcal{S}(\R^n) \Rightarrow f \phi \in \mathcal{S}(\R^n).$ 
\end{enumerate}
\end{lemma}

Heuristically, the functions $\phi \in \mathcal{S}(\R^n)$ \emph{decay rapidly at infinity} with all derivatives. Further, note that $\mathcal{S}(\R^n) \subset L^1(\R^n, \mathcal{B}_{\R^n}, d \bm{x})$, where $d\bm{x}$ is the Lebesgue measure on $\R^n$ \cite{Amir:master, slowly-growing-plus-prod-schwartz}. 
By the $L^1$ definition, each $\phi \in \mathcal{S}(\R^n)$ is measurable with respect to $(\R^n, \mathcal{B}_{\R^n})$ and absolutely integrable on $\R^n$; therefore, the necessary requirements for Fubini's theorem (see Thm.~1.4.28 in \cite{Amir:master}) are met:
\begin{equation}\label{Fubini}
\int_{\R^n} \phi(\x) \thinspace d \x = \int_{\R^{n-1}} \left[ \int_{\R} \phi(\bm{x}) \thinspace dx_j \right] d\bm{x}_{-j} \qquad \text{ for any } \phi \in \mathcal{S}(\R^n), \thinspace j \in \{1, \ldots, n\} 
\end{equation}

where $d\bm{x}_{-j}$ is the Lebesgue measure on $\R^{n-1}$ and $\int_{\R} \phi(\x) \thinspace d x_j$ is guaranteed to be finite. \emph{Note that this caution is warranted}, because Fubini's theorem does \emph{not} hold when any of the integrals involved are infinite; for one of several classic counter-examples, see \cite{fubini_counter_example}.


For a slowly growing $f \in \mathcal{P}(\R^n)$, Schwartz $\phi \in \mathcal{S}(\R^n)$, and $j \in [n],$ define the integral
\begin{equation}
I(\x) := \int_{\R^n} f(\x) \partial^{\bm{e}_j} (\phi(\x)) \thinspace d\x,
\end{equation}
where $\bm{e}_j$ is the $j$th unit multi-index (\cref{math_background}). 
By Lemma \ref{slowly-growing-lemma}, the integrand $f \partial^{\bm{e}_j} (\phi) \in \mathcal{S}(\R^n)$, so: (1) $I(\bm{x})$ is finite, and (2) Fubini's theorem (Equation~\ref{Fubini}) gives
\begin{equation}
I(\x) = \int_{\R^{n-1}} \left[ \int_{\R} f(\x) \partial^{\bm{e}_j} (\phi(\x))  \thinspace dx_j \right] d\bm{x}_{-j}. 
\end{equation} 	
Since $f \partial^{\bm{e}_j} \phi \in \mathcal{S}(\R^n)$, when projecting onto one dimension (i.e., holding all $x_i$ for all $i \neq j$ fixed), the integrand of $\int_{-R}^R  f(\x) \partial^{\bm{e}_j} (\phi(\x))  \thinspace dx_j$ belongs to $\mathcal{S}(\R)$ (where $R > 0$ is arbitrary). We can also note that Schwartz functions are Lipschitz continuous  \cite{MSE:derivSchwartzuniforcontinuous}
and Lipschitz functions on a compact interval $[-R, R]$ are absolutely continuous (Problem 5.20a in \cite{LipCompactAbsCont}) and thus uniformly continuous. Therefore, the needed conditions for integration by parts are met \cite{h_grad,conds_for_IBP}. By Integration by Parts in one dimension,

\begin{align*}
\int_{\R} f(\x) \partial^{\bm{e}_j} (\phi(\x)) &= \lim_{R \rightarrow \infty} \int_{-R}^R  f(\x) \partial^{\bm{e}_j} (\phi(\x))  \thinspace dx_j \\
&= \lim_{R \rightarrow \infty} \bigg[f(\bm{x})\phi(\bm{x}) \bigg]_{x_j = -R}^{x_j =+R} - \lim_{R \rightarrow \infty} \int_{-R}^R \partial^{e_j} f(\bm x) \phi(\bm x) \thinspace d x_j \\
&= - \int_{\R} \partial^{e_j} f(\bm x) \phi(\bm x) \thinspace d x_j.
\end{align*}


It is a straight-forward exercise to show that the boundary terms vanish. Hence, by applying Fubini (Eq.~\ref{Fubini}) again,
\begin{equation*}
I(\x) = \int_{\R^{n-1}} \left[- \int_{\R}  \partial^{e_j} f(\bm x) \phi(\bm x) \thinspace dx_j \right] d\bm{x}_{-j} = - \int_{\R^n} \partial^{e_j} f(\bm x) \phi(\bm x) \thinspace d\x
\end{equation*} 
and so
\begin{equation}\label{IBP_one_step}
\int_{\R^n} f(\x) \partial^{\bm{e}_j} (\phi(\x)) \thinspace d\x = - \int_{\R^n} \partial^{e_j} f(\bm x) \phi(\bm x) \thinspace d\x  \text{ for any } f \in \mathcal{P}(\R^n), \phi \in \mathcal{S}(\R^n), \thinspace j \in [n],
\end{equation}
with both integrals finite. 


The following statement is well-know, but we include it as a theorem to confirm that the integrals used in later proofs are well-defined (e.g., finite). 

\begin{theorem}\label{schwartz_gives_moments_exist} For $\bxi \sim \mathcal{N}(\bm{0}, \Sigma)$, $\mathbb{E}_{\bxi}(f(\bm{\xi}))$ exists for all $f \in \mathcal{P}(\R^n)$.
\end{theorem}
\begin{proof} By definition of existing expectation, we must prove that that $f(\bxi) \in L^1(\Omega, \mathcal{A}, \mathbb{P})$, which involves showing that: (1)  $f(\bxi)$ is $(\Omega, \mathcal{A})$-measurable, and (2)  $\int_\Omega |f(\bxi(\omega))| \thinspace d\mathbb{P}(\omega) < \infty.$ We of course have the map 
\[ (\Omega, \mathcal{A}, \mathbb{P}) \overset{\bm{\xi}}{\longrightarrow} (\R^n, \mathcal{B}_{\R^n}, d \bm{\xi}) \overset{f}{\longrightarrow} (\R, \mathcal{B}_{\R}, d \bm{x}) \] 
where $d\bm{\xi}$ is the \emph{law} of $\bm{\xi}$ and $d \bm{x}$ is the Lebesgue measure on $\R$. 
Since $\bm{\xi}$ is a random variable (and thus $(\Omega, \mathcal{A})$-measurable by definition), and $f$ is continuous, the composition $f(\bm{\xi})$ is $(\Omega, \mathcal{A})$-measurable as well (follows immediately since continuous functions have open preimages, and the Borel sigma algebra on $\R^n$ is generated by open balls, e.g., cite Proposition 1 in \cite{PETROVAI2019640}). By Corollary 1.3.62 in \cite{Amir:master},
$\int_{\Omega} |f(\bm{\xi} (\omega))| \thinspace d\mathbb{P}(\omega) < \infty$ if and only if $\int_{\R^n}  |f(\bm{x})| \thinspace\varphi(\bm{x})  d\bm{x} < \infty$, where $\varphi: \R^n \R_{\ge 0}]$ is the PDF of $\bxi$. Since the PDF $\varphi \in \mathcal{S}(\R^n)$ and $f \in \mathcal{P}(\R^n)$, by Lemma \ref{slowly-growing-lemma}, $f \varphi \in \mathcal{S}(\R^n) \subset L^1(\R^n, \mathcal{B}_{\R^n}, d \x)$; thus, $ \int_{\R^n}  |f(\bm{x})| \thinspace\varphi(\bm{x})  d\bm{x}  = \int_{\R^n} |f(\x) \varphi(\x)| \thinspace d\x  < \infty$, completing the proof. 
\end{proof}
\begin{corollary}\label{moments-exist} Let $0 \le N < \infty$ be an integer. The $N$th order triple product 
$$\langle H_{\balpha^{(1)}}, \ldots, H_{\balpha^{(N)}} \rangle = \mathbb{E}_{\bm{\xi}}\left(\prod_{i = 1}^N H_{\balpha^{(i)}}(\bm{\xi}) \right)$$ exists. 
\end{corollary}
\begin{proof}
Since each $H_{\balpha^{(i)}}:\R^n \rightarrow \R$ is a polynomial, clearly $H_{\balpha^{(i)}}:\R^n \rightarrow \R \in \mathcal{P}(\R^n)$; hence, by Lemma \ref{slowly-growing-lemma}, their product $\prod_{i = 1}^N H_{\balpha^{(i)}} \in \mathcal{P}(\R^n)$ as well. By Thm.~\ref{schwartz_gives_moments_exist}, $\mathbb{E}_{\bm{\xi}}\left( \prod_{i = 1}^N H_{\balpha^{(i)}}(\bm{\xi}) \right) $ exists.
\end{proof}

\section{Existence of higher-order products}\label{existence_section}
With these technical results available, we now demonstrate that the $N$th order product, along with the expected values involved in the proof of its formula (Thm.~\ref{N_prod}), are well-defined. That is, we show that the integrals used later in this document are finite and their integrands are measurable, so they belong to $L^1(\Omega, \mathcal{A}, d \bm{\xi}).$\footnote{One might be concerned that membership to $L^1(\Omega, \mathcal{A}, d \bm{\xi})$ is not restrictive enough, given that polynomial chaos expansions rely on finite variance. However, finite variance is the restriction on the \emph{stochastic process} being represented via the CPCE, not on the integrals involved in the the forthcoming proofs of the $N$th order products. As we will establish in Sect.~\ref{Nthprodformula}, the basis polynomials $H_{\balpha}(\bxi)$ themselves have all finite moments of all orders.} In order to do so, we need any functions of $\bm{\xi}$ and their derivatives to have at most polynomial growth, which allows their products with the PDF $\varphi$ of $\bxi$ to be Schwartz (Def.~\ref{Schwartz}) --- and thus decay fast enough for their integrals to be finite. We refer to smooth, real-valued functions whose derivatives grow at most polynomially as \emph{slowly growing} via Definition \ref{slowly_increasing}.

First, we first establish that the expected value with respect to the $d\bm{\xi}$ measure exists for any slowly growing statistic of the Gaussian $\bxi$ random variable (Theorem \ref{schwartz_gives_moments_exist}). As a corollary of this theorem (Corollary \ref{moments-exist}), we have that the $N$th order moments of the multivariate Hermite polynomials exist as well. Then, using Fubini's theorem/integration by parts (\cref{Fubini} in \cref{appendix_defs}), we can prove one of the central results (Theorem \ref{IBP}) of this paper, which states that the expected value of a MHP multiplied by a slowly growing function can be expressed as the expected value of a partial derivative of the slowly growing function; that is,
\begin{align}
\mathbb{E}_{\bxi}(H_{\A}(\bxi) f(\bxi)) = \mathbb{E}_{\bxi}(\partial^{\A} f(\bxi)) \qquad \text{ for $f$ slowly growing}.
\end{align}
Equipped with Theorem \ref{IBP}, we finally prove in Corollary \ref{general_deriv_N_product} that the $N$th order product is finite and can be written as 
\begin{align}\label{stuff}
\langle H_{\balpha^{(1)}}, \ldots, H_{\balpha^{(N)}} \rangle = \mathbb{E}_{\bxi}\left( \partial^{\balpha^{(j)}} \left( \textstyle\prod_{i \neq j} H_{\balpha^{(i)}}(\bxi) \right) \right).
\end{align}
Essentially, (\ref{stuff}) shows that we can ``transfer" any one of the MHPs $H_{\balpha^{(j)}}$ from the $N$th product (left-hand side) over to the partial derivative operator $\partial^{\balpha^{(j)}}$ in expectation (right-hand side).


\begin{theorem}\label{IBP}Let $f \in C^\infty(\R^n, \R)$ be slowly growing. Let $\A \in \N_0^n$ be a fixed but arbitrary multi-index. Then
\[ \mathbb{E}_{\bm \xi} (H_{\A}(\bm \xi) f(\bm \xi))  = \mathbb{E}_{\bm \xi}( \partial^{\A} f(\bm \xi)). \] 
\end{theorem}
\begin{proof} First note that $ \mathbb{E}_{\bm \xi} (H_{\A}(\bm \xi) f(\bm \xi))$ and $ \mathbb{E}_{\bm \xi} (H_{\A}(\bm \xi) f(\bm \xi))$ exist by Lemma \ref{slowly-growing-lemma} and Theorem \ref{schwartz_gives_moments_exist}. If $\A = \bm{0}$, then the result is immediate, since $H_{\bm{0}} = 1.$ So suppose $\A \neq \bm{0}$. Then $A_j > 0$ for some $j \in [n].$ Letting $\varphi$ be the density of $\bxi$, we have
\begin{align*}
\mathbb{E}_{\bxi} (H_{\A}(\bxi) f(\bxi))
&= \int_{\Omega} H_{\A}(\bxi(\omega)) f(\bxi(\omega)) \, d \mathbb{P}(\omega) \\
&= \int_{\R^n} H_{\A}(\x) f(\x) \varphi(\x) \, d\x \\
&= (-1)^{|\A|} \int_{\R^n} f(\x)\, \partial^{\A} \varphi(\x) \, d\x && \text{[$H_{\A} \varphi = (-1)^{|\A|} \partial^{\A} \varphi$]} \\
&= (-1)^{|\A|} \int_{\R^n} f(\x)\, \partial^{\bm{e}_j} \big(\partial^{\A - \bm{e}_j} \varphi(\x)\big) \, d\x \\
&= (-1)^{|\A|+1} \int_{\R^n} \partial^{\bm{e}_j} f(\x)\, \partial^{\A - \bm{e}_j} \varphi(\x) \, d\x && \text{[by \cref{IBP_one_step}]}.
\end{align*}

Since $f \in \mathcal{P}(\R^n)$ and $\partial^{\A - \bm{e}_j} \varphi \in \mathcal{S}(\R^n)$ is Schwartz (because the Gaussian density $\varphi$ is itself Schwartz), we could safely apply Eq.~\ref{IBP_one_step} --- a consequence of integration by parts --- in the last equality above; namely, the integrals involved are finite. In fact, we can apply Eq.~\ref{IBP_one_step} a total of $|\A|$ times, because $\partial^{\B}f \in \mathcal{P}(\R^n)$ for any $\B \in \R^n$ (Lemma \ref{slowly-growing-lemma}) and $\partial^{\B} \varphi \in \mathcal{S}(\R^n)$ for any $\B \in \N_0^n.$ Furthermore, note that each $\int_{\R^n} \partial^{\B} f(\x) \partial^{\A - \B} \varphi(\x) \thinspace d \x < \infty$ since the integrand is Schwartz. Thus, by applying Eq.~\ref{IBP_one_step} an additional $|\A| - 1$ times, we have
\begin{align*}
 \mathbb{E}_{\bxi} (H_{\A}(\bm \xi) f(\bm \xi)) &= (-1)^{|\A| + |\A|} \int_{\R^n} \partial^{\A} f(\x) \partial^{\bm{0}} \varphi(\x) \thinspace d \x \\
 &= \int_{\R^n} \partial^{\A} f(\x) \varphi(\x) \thinspace d\x \\
 &= \mathbb{E}_{\bxi} (\partial^{\A} f(\bxi))
 \end{align*}
as desired.
\end{proof}
Note that Theorem \ref{IBP} \emph{relies on $\bm{\xi}$ being Gaussian} in order to have Schwartz integrands and therefore guarantee that all the integrals involved are finite; this finiteness is an \emph{essential} requirement for applying Fubini's theorem.\footnote{Perhaps a similar argument could apply for non-Gaussian random variables with compact support to ensure finite integrals, though several other complications would likely need to be considered. This is an area for future inquiry.}

Given Theorem \ref{IBP}, Corollary \ref{general_deriv_N_product} immediately gives us a formula for the $N$th order product.

\begin{corollary}\label{general_deriv_N_product}By Cor.~\ref{moments-exist}, we know the $N$th order product is finite. In addition, we have 
\begin{equation}
\langle H_{\balpha^{(1)}}, \ldots, H_{\balpha^{(N)}} \rangle = \mathbb{E}_{\bxi}\left( \partial^{\balpha^{(j)}} \left( \textstyle\prod_{i \neq j} H_{\balpha^{(i)}}(\bxi) \right) \right) \label{eq_gen_N_deriv_prod}
\end{equation}
for any $j \in [n].$ 
\end{corollary}
\begin{proof} Define $f(\bxi) = \textstyle\prod_{i \neq j} H_{\balpha^{(i)}}(\bxi) \in \mathcal{P}(\R^n)$. Then
\[\langle H_{\balpha^{(1)}}, \ldots, H_{\balpha^{(N)}} \rangle = \mathbb{E}_{\bxi} (H_{\balpha^{(j)}}(\bxi) f(\bxi))\]
Apply Theorem \ref{IBP}, and the result follows. 
\end{proof}
Corollary \ref{general_deriv_N_product} is \emph{extremely} useful, because it allows us to use known results about $(N-1)$th products --- which will appear on the right-hand side of Eq.~\ref{eq_gen_N_deriv_prod} --- when computing $N$th products. We now have all the needed ingredients for an iterative formula for $N$th order product of multivariate Hermite polynomials in correlated Gaussian random variables.



\section*{Code} Code is available from the author upon request. 


\Urlmuskip=0mu plus 1mu\relax

\bibliographystyle{siamplain}
\bibliography{nthprodbib}

\begin{thebibliography}{10}

\bibitem{alexanderian2013:master}
{\sc A.~Alexanderian}, {\em On spectral methods for variance based sensitivity
  analysis}, Probability Surveys, 10 (2013), pp.~51 -- 68,
  \url{https://doi.org/10.1214/13-PS219}.

\bibitem{fubini_counter_example}
{\sc J.~Boardman}, {\em A {F}ubini counterexample}.
\newblock 2001, \url{https://math.jhu.edu/~jmb/note/nofub.pdf} (accessed
  2025-07-09).

\bibitem{exceptionsToPCefficiency1}
{\sc M.~Branicki and A.~J. Majda}, {\em Fundamental limitations of polynomial
  chaos for uncertainty quantification in systems with intermittent
  instabilities}, Communications in Mathematical Sciences, 11 (2013),
  pp.~55--103, \url{https://api.semanticscholar.org/CorpusID:16416936}.

\bibitem{cammartin:master}
{\sc R.~H. Cameron and W.~T. Martin}, {\em The orthogonal development of
  non-linear functionals in series of {F}ourier-{H}ermite functionals}, Annals
  of Mathematics, 48 (1947), pp.~385 -- 392,
  \url{http://www.jstor.org/stable/1969178}.

\bibitem{constantine_primer:master}
{\sc P.~Constantine}, {\em A primer on stochastic {G}alerkin methods}, Feb.
  2007,
  \url{https://www.cs.colorado.edu/~paco3637/docs/constantine2007primer.pdf}
  (accessed 2020-02-06).

\bibitem{constantine_thesis:master}
{\sc P.~Constantine}, {\em Spectral methods for parametrized matrix equations},
  PhD thesis, Stanford University, 2009.

\bibitem{Debusschere2017}
{\sc B.~Debusschere}, {\em Intrusive Polynomial Chaos Methods for Forward
  Uncertainty Propagation}, Springer International Publishing, Cham, 2017,
  pp.~617 -- 636, \url{https://doi.org/10.1007/978-3-319-12385-1_19}.

\bibitem{Amir:master}
{\sc A.~Dembo}, {\em Probability theory: {STAT}310/{MATH}230}, 2019,
  \url{https://web.stanford.edu/class/stats310a/lnotes.pdf} (accessed
  2025-07-10).

\bibitem{dunkl_xu_2014:master}
{\sc C.~F. Dunkl and Y.~Xu}, {\em Orthogonal Polynomials of Several Variables},
  Encyclopedia of Mathematics and its Applications, Cambridge University Press,
  2nd~ed., 2014, \url{https://doi.org/10.1017/CBO9781107786134}.

\bibitem{ernst:master}
{\sc O.~G. Ernst, A.~Mugler, H.~J. Starkloff, and E.~Ullmann}, {\em On the
  convergence of generalized polynomial chaos expansions}, ESAIM: M2AN, 46
  (2012), pp.~317 -- 339, \url{https://doi.org/10.1051/m2an/2011045}.

\bibitem{Feinberg}
{\sc J.~Feinberg, V.~G. Eck, and H.~P. Langtangen}, {\em Multivariate
  polynomial chaos expansions with dependent variables}, SIAM Journal on
  Scientific Computing, 40 (2018), pp.~A199 -- A223,
  \url{https://doi.org/10.1137/15M1020447}.

\bibitem{exceptionsToPCefficiency2}
{\sc T.~Gerstner and M.~Griebel}, {\em Numerical integration using sparse
  grids}, Numerical Algorithms, 18 (1998), pp.~209--232,
  \url{https://doi-org.stanford.idm.oclc.org/10.1023/A:1019129717644}.

\bibitem{UQHandbook_IntroUQ}
{\sc R.~Ghanem, D.~Higdon, and H.~Owhadi}, {\em Introduction to Uncertainty
  Quantification}, Springer International Publishing, Cham, 2017, pp.~3--6,
  \url{https://doi.org/10.1007/978-3-319-12385-1_1}.

\bibitem{ghanem_and_spanos:master}
{\sc R.~G. Ghanem and P.~D. Spanos}, {\em Stochastic Finite Elements: A
  Spectral Approach}, Springer-Verlag, Berlin, Heidelberg, 1991.

\bibitem{Golub}
{\sc G.~H. Golub and C.~F.~V. Loan}, {\em Matrix Computations}, The Johns
  Hopkins University Press, 4th~ed., 2013.

\bibitem{h_grad}
{\sc H.~Grad}, {\em Note on n-dimensional {H}ermite polynomials},
  Communications on Pure and Applied Mathematics, 2 (1949), pp.~325--330,
  \url{https://doi.org/10.1002/cpa.3160020402}.

\bibitem{MSE:derivSchwartzuniforcontinuous}
{\sc
  \href{https://math.stackexchange.com/questions/1057422/are-derivatives-of-schwartz-functions-uniformly-continuous&amp7D&amp7B}{MathStackExchange}
  (Answered~by
  User:~\href{https://math.stackexchange.com/users/9849/davide-giraudo}{D.~Giraudo})},
  {\em Are derivatives of {S}chwartz functions uniformly continuous?}
\newblock Mathematics Stack Exchange, 2014,
  \url{https://math.stackexchange.com/questions/1057422/are-derivatives-of-schwartz-functions-uniformly-continuous&amp7D&amp7B}
  (accessed 2025-07-10).

\bibitem{conds_for_IBP}
{\sc
  \href{https://math.stackexchange.com/questions/2419391/lebesgue-spaces-and-integration-by-parts}{MathStackExchange}
  (Answered by User:
  \href{https://math.stackexchange.com/users/121671/mittens}{Dr.~Oliver Diaz}
  of University of Texas~at Austin)}, {\em {L}ebesgue spaces and integration by
  parts}.
\newblock Mathematics Stack Exchange, 2020,
  \url{https://math.stackexchange.com/questions/2419391/lebesgue-spaces-and-integration-by-parts}
  (accessed 2025-07-10).

\bibitem{MSE:prod-schwartz-poly}
{\sc
  \href{https://math.stackexchange.com/questions/4303036/product-of-a-schwartz-function-and-a-function-with-polynomial-bonuded-derivative?}{MathStackExchange}
  (Answered~by
  User:~\href{https://math.stackexchange.com/users/983565/dreamar}{DreamAR})},
  {\em Product of a {S}chwartz function and a function with polynomial bonuded
  [sic] derivatives is {S}chwartz?}
\newblock Mathematics Stack Exchange, Nov 2021,
  \url{https://math.stackexchange.com/questions/4303036/product-of-a-schwartz-function-and-a-function-with-polynomial-bonuded-derivative?noredirect=1&amp;lq=1}
  (accessed 2025-07-10).

\bibitem{slowly-growing-plus-prod-schwartz}
{\sc J.~K. Hunter}, {\em Notes on partial differential equations, appendix
  5{A}, page 134}, June 2014,
  \url{https://www.math.ucdavis.edu/~hunter/m218a_09/ch5A.pdf} (accessed
  2025-07-10).

\bibitem{janson_1997:master}
{\sc S.~Janson}, {\em {G}aussian {H}ilbert Spaces}, Cambridge Tracts in
  Mathematics, Cambridge University Press, 1997,
  \url{https://doi.org/10.1017/CBO9780511526169}.

\bibitem{mydatataable}
{\sc L.~Lyman}, {\em Higher order products numeric examples},
  \url{https://docs.google.com/spreadsheets/d/1XoJo_k-YSO772KoqtytVZFK4kJqttJG13DivRCkQ5uY/edit?usp=sharing}
  (accessed 2025-07-10).

\bibitem{mytripprodpaper}
{\sc L.~Lyman}, {\em Triple products of correlated multivariate {H}ermite
  polynomials}, 2026, \url{https://arxiv.org/abs/2610.30204}.

\bibitem{mypaper:master}
{\sc L.~Lyman and G.~Iaccarino}, {\em A bluff-and-fix algorithm for polynomial
  chaos methods}, in Computational Science -- ICCS 2020, V.~V.
  Krzhizhanovskaya, G.~Z{\'a}vodszky, M.~H. Lees, J.~J. Dongarra, P.~M.~A.
  Sloot, S.~Brissos, and J.~Teixeira, eds., Cham, 2020, Springer International
  Publishing, pp.~742 -- 756,
  \url{https://doi.org/0.1007/978-3-030-50436-6_55}.

\bibitem{mydoubprodpaper}
{\sc L.~Lyman and G.~Iaccarino}, {\em Second order moments of multivariate
  {H}ermite polynomials in correlated random variables}, in Computational
  Science -- ICCS 2021, M.~Paszynski, D.~Kranzlm{\"u}ller, V.~V.
  Krzhizhanovskaya, J.~J. Dongarra, and P.~M.~A. Sloot, eds., Cham, 2021,
  Springer International Publishing, pp.~698--712,
  \url{https://doi.org/10.1007/978-3-030-77980-1_53}.

\bibitem{le_maitre:master}
{\sc O.~L. Ma\^{i}tre and O.~Knio}, {\em Spectral Methods for Uncertainty
  Quantification: With Applications to Computational Fluid Dynamics}, January
  2010, \url{https://doi.org/10.1007/978-90-481-3520-2}.

\bibitem{stoch_galer_2:master}
{\sc H.~Matthies and A.~Keese}, {\em Galerkin methods for linear and nonlinear
  elliptic stochastic partial differential equations}, Computer Methods in
  Applied Mechanics and Engineering, 194 (2005), pp.~1295 -- 1331,
  \url{https://doi.org/10.1016/j.cma.2004.05.027}.

\bibitem{mueller2024polynomialchaossurrogateconstruction}
{\sc J.~N. Mueller, K.~Sargsyan, C.~J. Daniels, and H.~N. Najm}, {\em
  Polynomial chaos surrogate construction for random fields with parametric
  uncertainty}, 2024, \url{https://arxiv.org/abs/2311.00553},
  \url{https://arxiv.org/abs/2311.00553}.

\bibitem{Nataf}
{\sc A.~Nataf}, {\em Determination des distributions de probabilite\'{e}s dont
  les marges sont donn\'{e}es}, C.~R.~Acad.~Sci.~Paris, 255 (1962), pp.~42 --
  43.

\bibitem{Navarro}
{\sc M.~Navarro, J.~Witteveen, and J.~Blom}, {\em Polynomial chaos expansion
  for general multivariate distributions with correlated variables}, 2014,
  \url{https://doi.org/10.48550/ARXIV.1406.5483}.

\bibitem{Gauss_mult_idx_ref:master}
{\sc S.~Noreddine and I.~Nourdin}, {\em On the {G}aussian approximation of
  vector-valued multiple integrals}, Journal of Multivariate Analysis, 102
  (2011), pp.~1008 -- 1017, \url{https://doi.org/10.1016/j.jmva.2011.02.001}.

\bibitem{verif_valid_ref_2}
{\sc W.~Oberkampf and T.~Trucano}, {\em Validation methodology in computational
  fluid dynamics}, AIAA 2000-2549,  (2000),
  \url{https://doi.org/10.2514/6.2000-2549}.

\bibitem{mcbook:master}
{\sc A.~B. Owen}, {\em Monte {C}arlo theory, methods and examples}, 2013,
  \url{https://artowen.su.domains/mc/} (accessed 2025-07-10).

\bibitem{Digression_Herm_Polys}
{\sc K.~Patarroyo}, {\em A digression on {H}ermite polynomials}, 2019,
  \url{https://doi.org/10.48550/ARXIV.1901.01648},
  \url{https://arxiv.org/abs/1901.01648}.

\bibitem{Perko_2016}
{\sc Z.~Perkó, S.~R. van~der Voort, et~al.}, {\em Fast and accurate
  sensitivity analysis of impt treatment plans using polynomial chaos
  expansion}, Physics in Medicine \& Biology, 61 (2016), p.~4646,
  \url{https://doi.org/10.1088/0031-9155/61/12/4646},
  \url{https://dx.doi.org/10.1088/0031-9155/61/12/4646}.

\bibitem{PETROVAI2019640}
{\sc D.~M. Petrovai}, {\em Some examples of non-measurable {L}ebesgue
  functions}, Procedia Manufacturing, 32 (2019), pp.~640 -- 642,
  \url{https://doi.org/https://doi.org/10.1016/j.promfg.2019.02.264},
  \url{https://www.sciencedirect.com/science/article/pii/S2351978919302999}.
\newblock 12th International Conference Interdisciplinarity in Engineering,
  INTER-ENG 2018, 4–5 October 2018, Tirgu Mures, Romania.

\bibitem{GI_book}
{\sc M.~P. Pettersson, G.~Iaccarino, and J.~Nordström}, {\em Polynomial Chaos
  Methods for Hyperbolic Partial Differential Equations: Numerical Techniques
  for Fluid Dynamics Problems in the Presence of Uncertainties}, Springer
  International, Switzerland, 2015,
  \url{https://doi.org/https://doi.org/10.1007/978-3-319-10714-1}.

\bibitem{hutch}
{\sc R.~Pulch}, {\em Methods based on polynomial chaos for quadratic delay
  differential equations with random parameters}, Proceedings in Applied
  Mathematics and Mechanics, 26 (2026), p.~e70083,
  \url{https://onlinelibrary.wiley.com/doi/abs/10.1002/pamm.70083}.

\bibitem{rahman:master}
{\sc S.~Rahman}, {\em Wiener--{H}ermite polynomial expansion for multivariate
  {G}aussian probability measures}, Journal of Mathematical Analysis and
  Applications, 454 (2017), pp.~303 -- 334,
  \url{https://doi.org/10.1016/j.jmaa.2017.04.062},
  \url{http://www.sciencedirect.com/science/article/pii/S0022247X17304250}.

\bibitem{RAHMAN2019124538}
{\sc S.~Rahman}, {\em Dimensionwise multivariate orthogonal polynomials in
  general probability spaces}, Applied Mathematics and Computation, 362 (2019),
  p.~124538, \url{https://doi.org/10.1016/j.amc.2019.06.052}.

\bibitem{ZDR}
{\sc Z.~Rosario}, {\em Precision margin : first-principles margins for aircraft
  design under uncertainty}, PhD thesis, Stanford University, 2020.

\bibitem{Rosenblatt}
{\sc M.~Rosenblatt}, {\em {Remarks on a Multivariate Transformation}}, The
  Annals of Mathematical Statistics, 23 (1952), pp.~470 -- 472,
  \url{https://doi.org/10.1214/aoms/1177729394},
  \url{https://doi.org/10.1214/aoms/1177729394}.

\bibitem{LipCompactAbsCont}
{\sc H.~L. Royden}, {\em Real Analysis}, Macmillan, New York, 3rd~ed., 1988.

\bibitem{Baby_Rudin}
{\sc W.~Rudin}, {\em Principles of Mathematical Analysis}, McGraw-Hill Book
  Company, Inc., New York-Toronto-London, 1953.

\bibitem{salazar}
{\sc J.~Salazar}, {\em The Polynomial Chaos Method with Applications to Random
  Differential Equations}, PhD thesis, The University of Texas at Arlington,
  2013.

\bibitem{proton_therapy_thesis}
{\sc J.~Salverda}, {\em Probabilistic treatment planning with polynomial chaos
  expansion for proton therapy}, master's thesis, Delft University of
  Technology, 2019.

\bibitem{leibniz-rule:master}
{\sc A.~Sasane}, {\em A Friendly Approach to Functional Analysis}, World
  Scientific Publishing Company, 2007, ch.~6.4, pp.~248 -- 249,
  \url{https://books.google.com/books/about/A_Friendly_Approach_To_Functional_Analys.html?id=vEQyDwAAQBAJ}.

\bibitem{slepian:master}
{\sc D.~Slepian}, {\em On the symmetrized {K}ronecker power of a matrix and
  extensions of {M}ehler’s formula for {H}ermite polynomials}, SIAM Journal
  on Mathematical Analysis, 3 (1972), pp.~606 -- 616,
  \url{https://doi.org/10.1137/0503060},
  \url{https://arxiv.org/abs/https://doi.org/10.1137/0503060}.

\bibitem{Smith}
{\sc R.~Smith}, {\em Uncertainty Quantification: Theory, Implementation, and
  Applications}, SIAM - Society for Industrial and Applied Mathematics, 2013.

\bibitem{airflowsoil}
{\sc P.~Sochala and O.~P. {Le Maître}}, {\em Polynomial chaos expansion for
  subsurface flows with uncertain soil parameters}, Advances in Water
  Resources, 62 (2013), pp.~139--154,
  \url{https://doi.org/https://doi.org/10.1016/j.advwatres.2013.10.003},
  \url{https://www.sciencedirect.com/science/article/pii/S030917081300184X}.

\bibitem{Soize_Ghanem}
{\sc C.~Soize and R.~Ghanem}, {\em Physical systems with random uncertainties:
  Chaos representations with arbitrary probability measure}, SIAM Journal on
  Scientific Computing, 26 (2004), pp.~395 -- 410,
  \url{https://doi.org/10.1137/S1064827503424505},
  \url{https://doi.org/10.1137/S106482750342450},
  \url{https://arxiv.org/abs/https://doi.org/10.1137/S1064827503424505}.

\bibitem{Stanley2001:master}
{\sc R.~P. Stanley}, {\em {Enumerative Combinatorics, Vol.~2}}, Cambridge
  University Press, 2001.

\bibitem{Lp_def:master}
{\sc E.~M. Stein and R.~Shakarchi}, {\em Functional Analysis: Introduction to
  Further Topics in Analysis}, Princeton University Press, 2012, ch.~1,
  \url{https://assets.press.princeton.edu/chapters/s9627.pdf}.

\bibitem{polys:master}
{\sc G.~Szeg{\H{o}}}, {\em Orthogonal Polynomials}, American Math. Soc:
  Colloquium publ, American Mathematical Society, 1975,
  \url{https://books.google.com/books?id=ZOhmnsXlcY0C}.

\bibitem{wiener:master}
{\sc N.~Wiener}, {\em The homogeneous chaos}, American Journal of Mathematics,
  60 (1938), pp.~897 -- 936, \url{http://www.jstor.org/stable/2371268}.

\bibitem{xiu_thesis:master}
{\sc D.~Xiu}, {\em Generalized (Wiener-Askey) Polynomial Chaos}, PhD thesis,
  Brown University, 2004.

\bibitem{xiu_stochcollsurvey}
{\sc D.~Xiu}, {\em Stochastic Collocation Methods: A Survey}, Springer
  International Publishing, Cham, 2017, pp.~699--716,
  \url{https://doi.org/10.1007/978-3-319-12385-1_26},
  \url{https://doi.org/10.1007/978-3-319-12385-1_26}.

\bibitem{xiu_karniadakis:master}
{\sc D.~Xiu and G.~E. Karniadakis}, {\em The {W}iener--{A}skey polynomial chaos
  for stochastic differential equations}, SIAM J. Sci. Comput., 24 (2002),
  pp.~619 -- 644, \url{https://doi.org/10.1137/S1064827501387826}.

\end{thebibliography}

\end{document}